\documentclass[a4paper,12pt]{article}
\usepackage{graphicx}
\usepackage{amsmath,amsthm}
\usepackage{amsfonts,amssymb}
\usepackage{amscd}
\usepackage{color}
\usepackage[cp1251]{inputenc}
\usepackage[matrix,arrow,curve]{xy}

\newtheorem{lemma}{Lemma}
\newtheorem{definition}{Definition}
\newtheorem{proposition}[lemma]{Proposition}
\newtheorem{theorem}[lemma]{Theorem}
\newtheorem{corollary}[lemma]{Corollary}
\newtheorem{Remark}{{\rm\bf Remark}}
\newcommand{\dimk}{{\rm dim_{\kk}}}
\newcommand{\Hom}{{\rm Hom}}
\newcommand{\Homk}{{\rm Hom}_{\kk}}
\newcommand{\charr}{{\rm char\,}}
\newcommand{\kk}{{\bf k}}
\newcommand{\dd}{{\bf d}}

\DeclareMathOperator{\Ker}{{Ker}}

\renewcommand{\Bar}{{\rm Bar}}

\newcommand{\Ext}{{\rm Ext}}

\DeclareMathOperator{\HH}{{HH}}
\DeclareMathOperator{\Ho}{{H}}
\renewcommand{\Im}{{\rm Im\,}}

\newcommand{\und}[1]{\overline{\text{\scriptsize $#1$}}}

\newcommand{\XX}{{\mathcal X}}

\newcommand{\BB}{{\mathcal B}}
\newcommand{\DD}{{\mathcal D}}

\newcommand{\DDD}{{\rm D}}
\newcommand{\KK}{{\rm K}}

\newcommand{\ii}{{\bf i}}
\newcommand{\ot}{\otimes}

\newcommand{\la}{\langle}
\newcommand{\ra}{\rangle}

\newcommand{\ee}{\epsilon}
\newcommand{\si}{\sigma}

\renewcommand{\le}{\leqslant}
\renewcommand{\ge}{\geqslant}

\newcommand{\ph}{\varphi}

\begin{document}
\title{Gerstenhaber bracket on the Hochschild cohomology via an arbitrary resolution}
\author{Yury Volkov}
%\email{wolf86\_666@list.ru}
%\classification{35J25}
%\keywords{Hochschild cohomology, Gerstenhaber bracket, bimodule resolution}
\date{}
\maketitle
\begin{abstract}
We prove formulas of different types that allow to calculate the Gerstenhaber bracket on the Hochschild cohomology of an algebra using some arbitrary projective bimodule resolution for it.
Using one of these formulas, we give a new short proof of the derived invariance of the Gerstenhaber algebra structure on Hochschild cohomology.
Also we give some new formulas for the Connes' differential on the Hochschild homology that lead to formulas for BV differential on the Hochschild cohomology in the case of symmetric algebras.
Finally, we use one of the obtained formulas to get a full description of the BV structure and, correspondingly, the Gerstenhaber algebra structure on the Hochschild cohomology of a class of symmetric algebras.
\end{abstract}
\section{Introduction}

Let $A$ be an associative unital algebra over a field $\kk$. The Hochschild
cohomology $\HH^*(A)$ of $A$ has a very rich structure. It  is a graded commutative algebra
via the cup product or the Yoneda product, and  it has a graded  Lie bracket  of degree $-1$ so that it becomes  a graded Lie
algebra; these make $\HH^*(A)$ a Gerstenhaber algebra \cite{Gerstenhaber}.
These structures have a good description in terms of the bar resolution of $A$, but this resolution is huge and so it is frequently useless for concrete computations.

The cup product is well studied. There are different formulas for computing it using an arbitrary projective resolution and they were used in many examples.
The situation with the Lie bracket is more complicated. Almost all computations of it are based on the method of so-called comparison morphisms. This method allows to transfer elements of Hochschild cohomology from one resolution to another.
For example, this method was applied for the description of the Lie bracket on the Hochschild cohomology of the group algebra of quaternion group of order 8 over a field of characteristic 2 in \cite{IIVZ}. Later this method was applied for all local algebras of the generalized quaternion type over a field of characteristic 2 in \cite{Ivan}. Applications of the method of comparison morphisms can be found also in \cite{ACT, RR, Sanch}.

Just a little time ago a formula for computing the bracket via a resolution which is not the bar resolution, appeared in \cite{NW}. The proof given there is valid for a resolution that satisfies some conditions. Other formulas for the Lie bracket are proved in the current work. These formulas use chain maps from a resolution to its tensor powers and homotopies for some null homotopic maps defined by cocycles. Then the formula of \cite{NW} is slightly changed and proved for an arbitrary resolution.
Note also that a nice formula for the bracket of a degree one element with an arbitrary element is given in \cite{SA}.

It is well known that the Hochschild cohomology is a derived invariant. The proof of this fact can be found, for example, in  \cite{Rickard}. The invariance of the cup product easily follows from this proof, while the derived invariance of the Gerstenhaber bracket was proved much later.
In \cite{Keller, Keller2} derived invariance of the Gerstenhaber bracket is proved using two different (relatively advanced) methods. In \cite{Keller} Keller employs the derived Picard group, while \cite{Keller2} relies on the use of DG categories. Here, using our new formulation of the bracket and the approach to the Hochschild homology proposed in \cite{Zimmermann}, we provide a direct proof of the derived invariance of the bracket which does not require any advanced technology.

Further, we give some formulas for the Lie bracket using so-called contracting homotopies. Then we discuss some formulas for the Connes' differential on the Hochschild homology. One of these formulas is a slight modification of the formula from \cite{Kaled}.
Also we give a formula using contracting homotopies for the Connes' differential. Thus, in the case where the Connes' differential induces a BV structure on Hochschild cohomology, we obtain an alternative way for the computing of the Lie bracket. We discuss this in the case where the algebra under consideration is symmetric.

Finally, we give an example of an application of the discussed formulas. We describe the BV structure and the Gerstenhaber bracket on the Hochschild cohomology of one family of symmetric local algebras of dihedral type. The Hochschild cohomology for these algebras was described in \cite{Gen1} and \cite{Gen2}. Note also that the Hochschild cohomology groups and the Hochschild cohomology ring modulo nilpotent radical were described in \cite{SnashTail} for a class of self-injective algebras including the family of symmetric algebras considered in this work. 

\section{Hochschild cohomology via the bar resolution}\label{HHB}

During this paper $A$ always denotes some algebra over a field $\kk$. We write simply $\otimes$ instead of $\otimes_{\kk}$.

 Let us recall how to define the Hochschild cohomology, the cup product and the Lie bracket in terms of the bar resolution. The Hochschild
cohomology groups are defined as $\HH^n(A)\cong \mathrm{Ext}^n_{A^e}(A, A)$ for $n \geq 0$, where $A^e=A \otimes
A^{\mathrm{op}}$ is the enveloping algebra of $A$.

\begin{definition}
An {\it $A^e$-complex} is a $\mathbb{Z}$-graded $A$-bimodule $P$ with a differential of degree $-1$, i.e. an $A$-bimodule $P$ 
with some fixed $A$-bimodule direct sum decomposition $P=\oplus_{n\in \mathbb{Z}}P_n$ and an $A$-bimodule homomorphism $d_P:P\rightarrow P$ such that $d_P(P_n)\subset 
P_{n-1}$ and $d_P^2=0$.  Let $d_{P,n}$ denote $d_P|_{P_{n}}$. The {\it $n$-th homology of $P$} is the vector space $\Ho_n(P)=(\Ker d_{P,n})/(\Im d_{P,n+1})$. 
An $A^e$-complex $P$ is called acyclic if $\Ho_n(P)=0$ for all $n\in\mathbb{Z}$ and is called {\it bounded on the right} if $P_n=0$ for small enough $n$. A {\it map of $A^e$-complexes} is a homomorphism of $A$-bimodules that 
respects the grading. If it also respects the differential, it is called a {\it chain map}.  A complex is called positive if $P_n=0$ for $n<0$. A pair $(P,\mu_P)$ is called a {\it resolution} of the algebra $A$ if $P$ is a positive complex, $\Ho_n(P)=0$ for $n>0$ and $\mu_P:P_0\to A$ is an $A$-bimodule homomorphism inducing an isomorphism $\Ho_0(P)\cong A$.
\end{definition}

Given an $A^e$-complex $P$, $(P,A)$ denotes the $\kk$-complex $\oplus_{n\le 0}\Hom_{A^e}(P_{-n},A)$ with differential $d_{(P,A),n}=\Hom_{A^e}(d_{P,-1-n},A)$.
Let $\mu_A: A\otimes A \to A$ be the multiplication map.

Let $\Bar(A)$ be the positive $A^e$-complex with $n$-th member $\Bar_n(A)=A^{\otimes(n+2)}$ for $n\ge 0$ and the differential $d_{\Bar(A)}$ defined by the equality
 $$
d_{\Bar(A)}(a_0\otimes\cdots\otimes a_{n+1})=\sum\limits_{i=0}^n(-1)^ia_0\otimes\cdots\otimes a_{i-1}\otimes a_ia_{i+1}\otimes 
a_{i+2}\otimes\cdots\otimes a_{n+1}
 $$
for $n>0$ and $a_i\in A$ ($0\le i\le n+1$). 
Then $(\Bar(A),\mu_A)$, is a projective 
$A^e$-resolution of $A$ that is called the bar resolution.

The Hochschild cohomology of the algebra $A$ is the homology of the complex $C(A)=(\Bar(A),A)$. 
We write $C^n(A)$ instead of $C_{-n}(A)$ and $\delta^n$ instead of $d_{C(A),-1-n}$.
Note that $C^0(A)\simeq A$ and $C^n(A)\simeq\Homk(A^{\otimes n},A)$.
Given $f\in C^n(A)$, we introduce the notation
$$
 \delta_n^i(f)(a_1\otimes\cdots\otimes a_{n+1}):=\begin{cases}
 a_1f(a_2\otimes\cdots\otimes a_{n+1}),&\mbox{if $i=0$},\\
  (-1)^if(a_1\otimes \cdots \otimes a_ia_{i+1} \otimes \cdots\otimes
a_{n+1}),&\mbox{if $1\le i\le n$},\\
(-1)^{n+1}f(a_1\otimes \cdots \otimes a_n) a_{n+1},&\mbox{if $i=n+1$}.
 \end{cases}
 $$
Then $\delta^n=\sum\limits_{i=0}^{n+1}\delta_n^i$. We have $\HH^n(A)=(\Ker \delta^n)/(\Im \delta^{n-1})$.

The cup product $\alpha\smile\beta \in
C^{n+m}(A)=\mathrm{Hom}_{\kk}(A^{\otimes (n+m)}, A)$ of $\alpha\in
C^n(A)$  and $\beta\in C^m(A)$ is given by
$$
(\alpha\smile\beta)(a_1\otimes \cdots\otimes
a_{n+m}):=\alpha(a_1\otimes\cdots\otimes a_n)
\beta(a_{n+1}\otimes\cdots\otimes a_{n+m}).
$$
This cup product induces a well-defined product in the Hochschild
cohomology
$$
\smile \colon \HH^n(A) \times \HH^m(A) \longrightarrow \HH^{n+m}(A)
$$
that turns the graded $\kk$-vector space $\HH^*(A)=\bigoplus_{n\geq
0}\HH^n(A)$ into a graded commutative algebra (\cite[Corollary
1]{Gerstenhaber}).

The Lie bracket is defined as follows. Let $\alpha \in C^n(A)$ and
$\beta \in C^m(A)$. If $n, m\geq 1$, then, for $1\leq i\leq n$, we define
$\alpha\circ_i \beta \in C^{n+m-1}(A)$ by the equality
$$
(\alpha\circ_i \beta)(a_1\otimes \cdots\otimes
a_{n+m-1}):=\alpha(a_1\otimes \cdots \otimes a_{i-1}\otimes
\beta(a_i\otimes \cdots \otimes a_{i+m-1})\otimes a_{i+m}\otimes
\cdots \otimes a_{n+m-1});
$$
if $ n\geq 1$ and $m=0$, then $\beta\in A$ and, for $1\leq i\leq
n$, we set
$$
(\alpha\circ_i\beta)(a_1\otimes \cdots\otimes
a_{n-1}):=\alpha(a_1\otimes \cdots \otimes a_{i-1}\otimes \beta
\otimes a_{i }\otimes \cdots \otimes a_{n-1});
$$
for any  other case, we set $\alpha\circ_i \beta$ to be zero. Now we
define
$$
\alpha\circ \beta
:=\sum_{i=1}^n(-1)^{(m-1)(i-1)}\alpha\circ_i\beta\mbox{ and }[\alpha,\, \beta] :=\alpha\circ \beta-(-1)^{(n-1)(m-1)}\beta\circ
\alpha.
$$
Note that $[\alpha,\, \beta]\in C^{n+m-1}(A)$. The operation $[\
\,,\,\ ]$ induces a well-defined Lie bracket on the Hochschild cohomology
$$
[\ \,,\,\ ] \colon \HH^n(A) \times \HH^m(A) \longrightarrow
\HH^{n+m-1}(A)
$$
such that $(\HH^*(A),\, \smile,\, [\ \,,\,\ ])$ is a Gerstenhaber
algebra (see \cite{Gerstenhaber}).

\section{Comparison morphisms}\label{CM}

Here we recall the method of comparison morphisms. But firstly we introduce some notation.

If $P$ is a complex, then we denote by $P[t]$ the complex, which equals to $P$ as an $A$-bimodule, with grading $P[t]_n=P_{t+n}$ and 
differential defined as $d_{P[t]}= (-1)^t d_{P}$. Note that $d_P$ defines a map from $P$ to $P[-1]$. Let now take some map of complexes $f:P\rightarrow Q$. For any $t\in\mathbb{Z}$, $f[t]$ denotes the map from $P[t]$ to $Q[t]$ induced by $f$, i.e. such a map that $f[t]|_{P[t]_i}=f|_{P_{i+t}}$. For simplicity we will write simply $f$ instead of $f[t]$, since in each situation $t$ can be easily recovered.
 Let $\dd f$ denote the map $fd_P-d_Qf:P\rightarrow Q[-1]$. We will frequently use the equality $\dd(fg)=(-1)^m(\dd f)g+f\dd g:N\rightarrow Q[m-1]$ that is valid for any $g:N\rightarrow P[m]$.
For two maps of complexes $f,g:P\rightarrow Q$ we write $f\sim g$ if $f-g=\dd s$ for some $s:U\rightarrow V[1]$. 
  Note that if $f\sim 0$ and $\dd g=0$, then $fg\sim 0$ and $gf\sim 0$ (for the composition that has sense). Also we always identify an $A$-bimodule $M$ with the complex $\tilde M$ such that $\tilde M_i=0$ for $i\not=0$ and $\tilde M_0=M$.
  Note also that if $f\sim 0$, then $\dd f=0$.
 It is not hard to see that if $P$ is a projective complex, $Q$ is exact in $Q_i$ for $i\ge n$, and $Q_i=0$ for $i<n$, then for any $f:P\rightarrow Q$ the equality $\dd f=0$ holds if and only if $f\sim 0$. Moreover, we have the following fact.
 
 \begin{lemma}\label{homlif}
 Let $P$ be a projective complex, $Q$ be exact in $Q_i$ for $i>n$, and $Q_i=0$ for $i<n$. Let  $\mu_Q:Q\rightarrow\Ho_n(Q)$ denote the canonical projection.
 If $f:P\rightarrow Q$ is such that $\dd f=0$ and $\mu_Qf\sim 0$, then $f\sim 0$.
 \end{lemma}
 \begin{proof}
 Assume that $\mu_Qf= \phi d_P$. Since $P_{n-1}$ is projective, there is some $\psi:P_{n-1}\rightarrow Q_n$ such that $\mu_Q\psi=\phi$. Then $f-\dd\psi$ is a chain map such that $\mu_Q(f-\dd\psi)=0$. Then it is easy to see that $f\sim\dd\psi\sim 0$.
 \end{proof}
 
 Let now $(P,\mu_P)$ and $(Q,\mu_Q)$ be two $A^e$-projective resolutions of $A$. The method of comparison morphisms is based on the following idea. Since $P$ is positive projective and $Q$ is exact in $Q_i$ for $i>0$, there is some chain map of complexes $\Phi_P^Q:P\rightarrow Q$ such that $\mu_Q\Phi_P^Q=\mu_P$. Analogously there is a chain map $\Phi_Q^P:Q\rightarrow P$ such that $\mu_P\Phi_Q^P=\mu_Q$. Then $\Phi_P^Q$ and $\Phi_Q^P$ induce maps from $(Q,A)$ to $(P,A)$ and backwards. Thus, we also have the maps 
 $$(\Phi_P^Q)^*:\Ho_*(Q,A)\rightarrow \Ho_*(P,A)\mbox{ and }(\Phi_Q^P)^*:\Ho_*(P,A)\rightarrow \Ho_*(Q,A).$$ Since $\dd\left(1_P-\Phi_Q^P\Phi_P^Q\right)=0$, we have $1_P\sim \Phi_Q^P\Phi_P^Q$ by the arguments above. Then it is easy to see that $(\Phi_P^Q)^*(\Phi_Q^P)^*=(\Phi_Q^P\Phi_P^Q)^*=1_{\Ho_*(P,A)}$ and, analogously, $(\Phi_Q^P)^*(\Phi_P^Q)^*=1_{\Ho_*(Q,A)}$. So we can define the Hochschild cohomology of $A$ as the homology of $(P,A)$, and this definition does not depend on the $A^e$-projective resolution $(P,\mu_P)$ of $A$. If we define some bilinear operation $*$ on $(Q,A)$, which induces an operation on $\HH^*(A)$, then we can define the operation $*_{\Phi}$ on $(P,A)$ by the formula $f*_{\Phi}g=(f\Phi_Q^P*g\Phi_Q^P)\Phi_P^Q$ for $f,g\in(P,A)$. It is easy to see that $*_{\Phi}$ induces an operation on $\HH^*(A)$ and that the induced operation coincides with $*$. Now we can take $Q=\Bar(A)$ and define the cup product and the Lie bracket on $(P,A)$ by the equalities
 $$
 f\smile_{\Phi}g=(f\Phi_{\Bar(A)}^P\smile g\Phi_{\Bar(A)}^P)\Phi_P^{\Bar(A)}\mbox{ and } [f,g]_{\Phi}=[f\Phi_{\Bar(A)}^P,g\Phi_{\Bar(A)}^P]\Phi_P^{\Bar(A)}.
 $$
 
 Thus, to apply the method of comparison morphism one has to describe the maps $\Phi_P^{\Bar(A)}$ and $\Phi_{\Bar(A)}^P$ and then use them to describe the bracket in terms of the resolution $P$. The problem is that for some $x\in P$ the formula $\Phi_P^{\Bar(A)}(x)$ is complicated and that to describe $\Phi_{\Bar(A)}^P$ one has to define it on a lot of elements.
 
Let now recall one formula for the cup product that uses an arbitrary $A^e$-projective resolution of  $A$ instead of the bar resolution. But firstly let us introduce some definitions and notation.

\begin{definition} Given $A^e$-complexes $P$ and $Q$, we define the {\it tensor product complex} $P\otimes_A Q$ by the equality $(P\otimes_A
Q)_n=\sum_{i+j=n}P_i\otimes_A Q_j$.  The differential $d_{P\otimes_A Q}$ is defined by the equality $d_{P\otimes_A Q}(x\otimes y)=d_P(x)\otimes
y+(-1)^ix\otimes d_Q(y)$ for $x\in P_i$, $y\in Q_j$.
\end{definition}

We always identify $P\ot_A A$ and $A\ot_A P$ with $P$ by the obvious isomorphisms of complexes. For any $n\in\mathbb{Z}$ we also identify $P\otimes_AQ[n]$ and $P[n]\otimes_AQ$ with $(P\otimes_AQ)[n]$. Note that this identification uses isomorphisms $\alpha_{P,Q}^n:P\otimes_AQ[n]\rightarrow (P\ot_A Q)[n]$ and $\beta_{P,Q}^n:P[n]\otimes_AQ\rightarrow (P\ot_A Q)[n]$ defined by the equalities
$\alpha_{P,Q}^n(x\ot y)=(-1)^{in}x\ot y$ and $\beta_{P,Q}(x\ot y)=x\ot y$ for $x\in P_i$ and $y\in Q$. In particular, we have two different isomorphisms $\beta_{P,Q}^n\alpha_{P[n],Q}^m$ and $\alpha_{P,Q}^m\beta_{P,Q[m]}^n$ from $P[n]\ot_AQ[m]$ to $(P\ot_A Q)[n+m]$. For convenience, we always identify $P[n]\ot_AQ[m]$ with $(P\ot_A Q)[n+m]$ using the isomorphism $\beta_{P,Q}^n\alpha_{P[n],Q}^m$ that sends $x\ot y$ to $(-1)^{(i+n)m}x\ot y$ for $x\in P_i$ and $y\in Q$. In particular, we identify $A[n]\ot_AA[m]$ to $A[n+m]$ by the isomorphism $\beta_{A,A}^n\alpha_{A[n],A}^m$  that sends $a\ot b$ to $(-1)^{mn}ab$ for $a,b\in A$.
 
\begin{definition}
Given an $A^e$-projective resolution $(P,\mu_P)$ of $A$, a chain map  $\Delta_P:P\rightarrow P^{\otimes_An}$ is called a {\it diagonal $n$-approximation} of $P$ if $\mu_P^{\otimes n}\Delta_P=\mu_P$.
\end{definition}

Let $(P,\mu_P)$ be an $A^e$-projective resolution of $A$.  Suppose also that $\Delta_P:P\rightarrow P\otimes_AP$ is a diagonal $2$-approximation of $P$. Then the operation $\smile_{\Delta_P}$ on $(P,A)$ defined for $f:P\rightarrow A[-n]$ and $g:P\rightarrow A[-m]$ by the equality $f\smile_{\Delta_P}g=(-1)^{mn}(f\otimes g)\Delta_P$ induces the cup product on $\HH^*(A)$. Note also that if $f\in C^n(A)$ and $g\in C^m(A)$, then the equality $f\smile g=(-1)^{mn}(f\otimes g)\Delta$ holds for $\Delta$ defined by the equality
\begin{equation}\label{Delta}
\Delta(1\otimes a_1\otimes\cdots\otimes a_n\otimes 1)=\sum\limits_{i=0}^n(1\otimes a_1\otimes\cdots a_i\otimes 1)\otimes_A(1\otimes a_{i+1}\otimes\cdots\otimes a_n\otimes 1).
\end{equation}

\section{Gerstenhaber bracket via an arbitrary resolution}\label{IT}

In this section we prove some new formulas for the Gerstenhaber bracket. The existence of these formulas is based on the following lemma.

\begin{lemma}\label{homot}
Let $(P,\mu_P)$ be an $A^e$-projective resolution of $A$ and $f:P\rightarrow A[-n]$ be such that $fd_P=0$. Then
$f\otimes 1_P-1_P\otimes f:P\otimes_A P\rightarrow P[-n]$ is homotopic to 0.
\end{lemma}
\begin{proof} It is easy to check that $\dd(f\otimes 1_P-1_P\otimes f)=0$. Since $\mu_P(\mu_P\otimes 1_P-1_P\otimes\mu_P)=0$, there is some map $\phi:P\otimes_A P\rightarrow P[1]$ such that $\mu_P\otimes 1_P-1_P\otimes\mu_P=\dd\phi$.
Then $\mu_P(f\otimes 1_P-1_P\otimes f)=-f\dd\phi\sim 0$ and so $f\otimes 1_P-1_P\otimes f\sim 0$ by Lemma \ref{homlif}.
\end{proof}

\begin{corollary}\label{lift_1}
Let $P$, $f$ be as above and $\Delta_P$ be some diagonal $2$-approximation of $P$. Then $(f\otimes 1_P-1_P\otimes f)\Delta_P:P\rightarrow P[-n]$ is homotopic to 0.
\end{corollary}
\begin{proof} Since $\dd\Delta_P=0$, everything follows directly from Lemma \ref{homot}.
\end{proof}

\begin{definition}
Let $P$, $f$ and $\Delta_P$ be as above. We call $\phi_f:P\rightarrow P[1-n]$ a {\it homotopy lifting} of $(f,\Delta_P)$ if $\dd\phi_f=(f\otimes 1_P-1_P\otimes f)\Delta_P$ and $\mu_P\phi_f+f\phi\sim 0$ for some $\phi:P\rightarrow P[1]$ such that $\dd\phi=(\mu_P\otimes 1_P-1_P\otimes\mu_P)\Delta_P$.
%If $\phi_f:P\otimes_AP\rightarrow P$ is such that $\phi_f^d=f\otimes 1_P-1_P^n\otimes f$ and $\mu_P\phi_f+f\phi=_d0$ for some homotopy $\phi$ for $\mu_P\otimes 1_P-1_P\otimes\mu_P$, then we say that $\phi_f$ is a homotopy lifting for $f$.
\end{definition}

One can show following the proofs of Lemmas \ref{homlif} and \ref{homot} that some homotopy lifting exists for any cocycle. Alternatively, the existence of some $\tilde\phi_f$ such that $\dd\tilde\phi_f=(f\otimes 1_P-1_P\otimes f)\Delta_P$ follows from Corollary \ref{lift_1} and, in particular there is some $\phi$ satisfying the equality from the definition of a homotopy lifting. Easy calculation shows that $\mu_P\tilde\phi_f+f\phi$ is a cocycle. Then there is $u:P\rightarrow P[1-n]$ such that $\dd u=0$ and $\mu_Pu=\mu_P\tilde\phi_f+f\phi$, and hence $\tilde\phi_f-u$ is a homotopy lifting.
Now we are ready to prove our first formula.

\begin{theorem}\label{main}
Let $(P,\mu_P)$ be an $A^e$-projective resolution of $A$ and $\Delta_P:P\rightarrow P\otimes_AP$ be a diagonal $2$-approximation of $P$.
Let $f:P\rightarrow A[-n]$ and $g:P\rightarrow A[-m]$ represent some cocycles.
Suppose that $\phi_f$ and $\phi_g$ are homotopy liftings for $(f,\Delta_P)$ and $(g,\Delta_P)$ respectively. Then the Gerstenhaber bracket of the classes of $f$ and $g$ can be represented by the class of the element
\begin{equation}
[f,g]_{\phi,\Delta}=(-1)^mf\phi_g+(-1)^{m(n-1)}g\phi_f.
\end{equation}
\end{theorem}
\begin{proof} We will prove the assertion of the theorem in three steps.

1. Let us prove that the operation induced on the Hochschild cohomology by $[,]_{\phi,\Delta_P}$  does not depend on the choice of $\Delta_P$ and $\phi$. We do this in two steps:\\
\begin{itemize}
\item If $\phi_g$ and $\phi_g'$ are two homotopy liftings for $g$, then $\dd(\phi_g-\phi_g')=0$ and $\mu_P(\phi_g-\phi_g')\sim g\ee$ for some chain map $\ee:P\rightarrow P[1]$. Then $\ee\sim 0$ and $\mu_P(\phi_g-\phi_g')\sim 0$. 
Hence, $\phi_g-\phi_g'\sim 0$ and $f\phi'_g\sim f\phi_g$. Analogously, $g\phi'_f\sim g\phi_f$ and so $[f,g]_{\phi',\Delta}\sim [f,g]_{\phi,\Delta}$.
\item Let $\Delta'_P$ and $\Delta_P$ be two diagonal $2$-approximations of $P$ and $\phi_f$ and $\phi_g$ be homotopy liftings for $(f,\Delta_P)$ and $(g,\Delta_P)$ correspondingly. Then $\Delta_P'=\Delta_P+\dd u$ for some $u$. Note that if $\dd\phi=(\mu_P\otimes 1_P-1_P\otimes\mu_P)\Delta_P$, then $\dd\big(\phi+(\mu_P\otimes 1_P-1_P\otimes\mu_P)u\big)=(\mu_P\otimes 1_P-1_P\otimes\mu_P)\Delta_P'$, and hence $\phi_f'=\phi_f+(f\otimes 1_P-1_P\otimes f)u$
and $\phi_g'=\phi_g+(g\otimes 1_P-1_P\otimes g)u$ are homotopy liftings for $(f,\Delta_P')$ and $(g,\Delta_P')$. Now we have
\begin{multline*}
[f,g]_{\phi',\Delta'}-[f,g]_{\phi,\Delta}=(-1)^mf(g\otimes 1_P-1_P\otimes g)u+(-1)^{m(n-1)}g(f\otimes 1_P-1_P\otimes f)u\\
=\big((-1)^{m+mn}g\otimes f-(-1)^{m}f\otimes g+(-1)^{m(n-1)+mn}f\otimes g-(-1)^{m(n-1)}g\otimes f\big)u=0.
\end{multline*}
\end{itemize}

2. Let us prove that the operation induced on the Hochschild cohomology does not depend on the choice of an $A^e$-projective resolution of $A$. Let $(Q,\mu_Q)$ be another $A^e$-projective resolution of $A$. Let $\Phi_P^Q:P\rightarrow Q$ and $\Phi_Q^P:Q\rightarrow P$ be comparison morphisms, and $\phi_{f\Phi_Q^P\Phi_P^Q}$ and $\phi_{g\Phi_Q^P\Phi_P^Q}$ be homotopy liftings for $(f\Phi_Q^P\Phi_P^Q,\Delta_P)$ and $(g\Phi_Q^P\Phi_P^Q,\Delta_P)$ correspondingly. It is not difficult to check that
$\phi_{f\Phi_Q^P}=\Phi_P^Q\phi_{f\Phi_Q^P\Phi_P^Q}\Phi_Q^P$ and $\phi_{g\Phi_Q^P}=\Phi_P^Q\phi_{g\Phi_Q^P\Phi_P^Q}\Phi_Q^P$ are homotopy liftings for $(f\Phi_Q^P,\Delta_Q)$ and $(g\Phi_Q^P,\Delta_Q)$ correspondingly in this case. Here $\Delta_Q$ denotes the map $(\Phi_P^Q\otimes\Phi_P^Q)\Delta_P\Phi_Q^P$.
Then
\begin{multline*}
[f\Phi_Q^P,g\Phi_Q^P]_{\phi,\Delta}=(-1)^mf\Phi_Q^P\Phi_P^Q\phi_{g\Phi_Q^P\Phi_P^Q}\Phi_Q^P+(-1)^{m(n-1)}g\Phi_Q^P\Phi_P^Q\phi_{f\Phi_Q^P\Phi_P^Q}\Phi_Q^P\\
=[f\Phi_Q^P\Phi_P^Q,g\Phi_Q^P\Phi_P^Q]_{\phi,\Delta_P}\Phi_Q^P=[f,g]_{\phi,\Delta}\Phi_Q^P.
\end{multline*}

3. Suppose now that $(P,\mu_P)=(\Bar(A),\mu_A)$ and $\Delta_P=\Delta$, where $\Delta$ is the map from \eqref{Delta}. Let us define
\begin{multline*}
\phi_g(1\otimes a_1\otimes\cdots a_{i+m-1}\otimes 1)\\=\sum_{j=1}^i(-1)^{(m-1)j-1}\otimes a_1\otimes \cdots \otimes a_{j-1}\otimes g(a_j\otimes \cdots \otimes a_{j+m-1})\otimes a_{j+m}\otimes\cdots \otimes a_{i+m-1}\otimes 1
\end{multline*}
and analogously for $\phi_f$. Then we have $(-1)^mf\phi_g+(-1)^{m(n-1)}g\phi_f=[f,g]$ by definition. Direct calculations show that $\phi_f$ and $\phi_g$ are homotopy liftings for $(f,\Delta)$ and $(g,\Delta)$ (in fact, $\phi_g$ coincides with $(-1)^mG(1_B\otimes g\otimes 1_B)\Delta^{(2)}$ in \cite[Notation 2.3]{NW} and the fact that $\phi_g$ is a homotopy lifting follows from \cite[Proposition 2.4]{NW} and our discussion below).
\end{proof}

Let $(P,\mu_P)$ be an $A^e$-projective resolution for $A$, and $\Delta_P^{(2)}:P\rightarrow P\otimes_AP\otimes_AP$ be a diagonal $3$-approximation of $P$. There is some homotopy $\phi_P$ for $\mu_P\otimes 1_P-1_P\otimes\mu_P$. Since 
$$(\mu_P\otimes\mu_P)(\mu_P\otimes 1_P\otimes 1_P-1_P\otimes 1_P\otimes \mu_P)\Delta_P^{(2)}=0,$$
there is some homotopy $\ee_P$ for $(\mu_P\otimes 1_P\otimes 1_P-1_P\otimes 1_P\otimes \mu_P)\Delta_P^{(2)}$. We define
\begin{equation}\label{brack}
f\circ_{\Delta_P^{(2)},\phi_P,\ee_P}g=f\phi_P(1_P\otimes g\otimes 1_P)\Delta_P^{(2)}-(-1)^{m}(f\otimes g)\ee_P:P\rightarrow A[1-n-m]
\end{equation}
and
$$
[f,g]_{\Delta_P^{(2)},\phi_P,\ee_P}=f\circ_{\Delta_P^{(2)},\phi_P,\ee_P}g-(-1)^{(n-1)(m-1)}g\circ_{\Delta_P^{(2)},\phi_P,\ee_P}f.
$$
This formula is a slightly corrected variant of the formula from \cite{NW}.

\begin{corollary}\label{NWF}
The operation $[,]_{\Delta_P^{(2)},\phi_P,\ee_P}$ induces an operation on $\HH^*(A)$ that coincides with the usual Lie bracket on the Hochschild cohomology.
\end{corollary}
\begin{proof}
By Theorem \ref{main} it is enough to check that $-(1_P\otimes g)\ee_P+(-1)^m\phi_P(1_P\otimes g\otimes 1_P)\Delta_P^{(2)}$ is a homotopy lifting for $(g,(\mu_P\otimes 1_P\otimes 1_P)\Delta_P^{(2)})$ if $gd_P=0$. Let us verify the first condition:
\begin{multline*}
-\dd\big((1_P\otimes g)\ee_P+(-1)^m\phi_P(1_P\otimes g\otimes 1_P)\Delta_P^{(2)}\big)=-(1_P\otimes g)\dd\ee_P+\dd\phi_P(1_P\otimes g\otimes 1_P)\Delta_P^{(2)}\\
=(1_P\otimes g)(1_P\otimes 1_P\otimes \mu_P-\mu_P\otimes 1_P\otimes 1_P)\Delta_P^{(2)}+(\mu_P\otimes 1_P-1_P\otimes\mu_P)(1_P\otimes g\otimes 1_P)\Delta_P^{(2)}\\
=(g\ot 1_P-1_P\otimes g)(\mu_P\otimes 1_P\otimes 1_P)\Delta_P^{(2)}.
\end{multline*}
The second condition can be easily verified after noting that $\Im\phi_P\subset\oplus_{i>0}P_i\subset\Ker\mu_P$. Indeed, we have
\begin{multline*}
\mu_P\big(-(1_P\otimes g)\ee_P+(-1)^m\phi_P(1_P\otimes g\otimes 1_P)\Delta_P^{(2)}\big)+g\phi_P(\mu_P\otimes 1_P\otimes 1_P)\Delta_P^{(2)}\\
=g\big(\phi_P(\mu_P\otimes 1_P\otimes 1_P)\Delta_P^{(2)}-(\mu_P\ot 1_P)\ee_P\big)\sim 0
\end{multline*}
because $\dd\big(\phi_P(\mu_P\otimes 1_P\otimes 1_P)\Delta_P^{(2)}-(\mu_P\ot 1_P)\ee_P\big)=0$.
\end{proof}

\begin{Remark}{\rm
Usually the diagonal $3$-approximation $\Delta_P^{(2)}$ is constructed using some $2$-approximation $\Delta_P$ by the rule $\Delta^{(2)}_P=(\Delta_P\otimes 1_P)\Delta_P$. It often occurs that the maps $\Delta_P$ and $\mu_P$ satisfy the equality 
\begin{equation}\label{goodcond}
(\mu_P\otimes 1_P)\Delta_P=1_P=(1_P\otimes \mu_P)\Delta_P.
\end{equation}
In this case some things become easier. Firstly, one can set $\phi=0$ in the definition of a homotopy lifting. Then the second condition simply means that $\mu_P\phi_f$ is a coboundary. In particular, one can simply set $\phi_f|_{P_{n-1}}=0$. Secondly, if \eqref{goodcond} holds and the diagonal $3$-approximation is defined as above, then one can set $\ee_P=0$ in equality \eqref{brack}. Thus, we get the formula from \cite{NW} in the case where \eqref{goodcond} holds. Note that the condition \eqref{goodcond} is weaker than the conditions proposed in \cite{NW}.

On the other hand, we always can set $\ee_P=(\phi_P\otimes 1_P+1_P\otimes\phi_P)\Delta_P^{(2)}$ and obtain the following formula for the bracket:
\begin{multline}\label{brack2}
[f,g]_{\Delta_P^{(2)},\phi_P,\ee_P}=-f\phi_P(g\otimes 1_P\otimes 1_P-1_P\otimes g\otimes 1_P+1_P\otimes 1_P\otimes g)\Delta_P^{(2)}\\
+(-1)^{(n-1)(m-1)}g\phi_P(f\otimes 1_P\otimes 1_P-1_P\otimes f\otimes 1_P+1_P\otimes 1_P\otimes f)\Delta_P^{(2)}.
\end{multline}
}
\end{Remark}

\begin{Remark}{\rm
In fact, Corollary \ref{lift_1} can be proved directly without Lemma \ref{homot}. Then one can show that homotopy liftings exist using only the projectivity  of $P$ and not of its tensor powers over $A$. This allows to define the Gerstenhaber bracket on $\mathrm{Ext}^*_{A^e}(A, A)$ for any associative ring $A$ even in the case where one cannot use the bar resolution for this.
}
\end{Remark}

%\begin{Remark}{\rm
%If $\Phi_1,\Phi_2:P\rightarrow P$ are two comparison morphisms and $\phi_P$ and $\ee_P$ satisfy the conditions $\mu_P\otimes \Phi_1-\Phi_2\otimes\mu_P=\phi_P^d$ and $(\mu_P\otimes 1_P \otimes \Phi_1-\Phi_2\otimes 1_P\otimes \mu_P)\Delta_P=\ee_P^d$, then one can %show that $-(1_P^m\otimes g)\ee_P+(-1)^m\phi_P(1_P^m\otimes g\otimes 1_P)\Delta_P$ is a homotopy lifting for $(g,(\mu_P\otimes 1_P\otimes \Phi_1)\Delta_P)$ if $gd_P=0$. So the statement of the Corollary \ref{NWF} holds for such $\phi_P$ and $\ee_P$.
%}
%\end{Remark}

\section{Derived invariance of the Gerstenhaber bracket}

Let $\DDD^- A$ and $\KK_p^-A$ denote the derived category of bounded on the right complexes of $A$-modules and the homotopy category of bounded on the right complexes of  $A$-projective modules respectively. Note that the construction of a projective resolution for a complex induces an equivalence between $\DDD^- A$ and $\KK_p^-A$.
In this section $(P,\mu_P)$ is called a projective bimodule resolution of $A$ if $P\in\KK_p^-A^e$ and the morphism of $A^e$-complexes $\mu_P:P\rightarrow A$ induces an isomorphism in homology, i.e. $P$ does not have to be concentrated only in nonnegative degrees.
Then the chain map $\Delta_P:P\rightarrow P\ot_AP$ is called a diagonal $2$-approximation of $P$ if $(\mu_P\ot_A\mu_P)\Delta_P\sim\mu_P$.

One can easily check that all the arguments of the previous sections are valid in the settings of this section. In particular, for any map $f:P\rightarrow A[-n]$ there exists a homotopy lifting for $(f,\Delta_P)$ and the statement of Theorem \ref{main} holds.

We will say that $A$ is {\it standardly derived equivalent} to $B$ if there exist $U\in \DDD^-(A\ot B^{\rm op})$ and $V\in\DDD^-(B\ot A^{\rm op})$ such that $U\ot_B^LV\cong A$ in $\DDD^-A^e$ and $V\ot_{A}^LU\cong B$ in $\DDD^-B^e$. We will assume without loss of generality that $U\in \KK_p^-(A\ot B^{\rm op})$ and $V\in\KK_p^-(B\ot A^{\rm op})$. The paper \cite{Rickard} guarantees that if $A$ and $B$ are algebras over a field, then they are standardly derived equivalent if and only if they are derived equivalent.
Since $U\in \KK_p^-(A\ot B^{\rm op})$, $V\in\KK_p^-(B\ot A^{\rm op})$, $U\ot_{B}^LV\cong A$ in $\DDD^-A^e$ and $V\ot_{A}^LU\cong B$ in $\DDD^-B^e$, there are chain maps $\alpha:U\ot_{B}V\rightarrow A$ and $\beta:V\ot_{A}U\rightarrow B$ that induce isomorphisms in homology. 
We will need the following technical lemmas.

\begin{lemma}\label{comp}
The maps $\alpha$ and $\beta$ above can be chosen in such a way that
$$\alpha\ot 1_U\sim 1_U\ot\beta:U\ot_B V\ot_A U\rightarrow U\mbox{ and }1_V\ot\alpha\sim \beta\ot 1_V:V\ot_AU\ot_{B}V\rightarrow V.$$
\end{lemma}
\begin{proof} Let $\tilde\beta:V\ot_{A}U\rightarrow B$ be some chain map inducing isomorphism in homology. Note that $\alpha(\alpha\ot 1_{U\ot_BV}-1_{U\ot_BV}\ot\alpha)=0$. Since $\alpha$ is a quasi-isomorphism, we have
$$\alpha\ot 1_{U\ot_BV}\sim 1_{U\ot_BV}\ot\alpha:U\ot_B V\ot_A U\ot_B V\rightarrow U\ot_B V.$$
Analogously, $\tilde\beta\ot 1_{V\ot_AU}\sim 1_{V\ot_AU}\ot \tilde\beta$. Let $\beta$ be a chain map that equals
$$\tilde\beta(1_V\ot\alpha\ot 1_U)(1_{V\ot_AU}\ot \tilde\beta^{-1})$$
in $\Hom_{\DDD^-B^e}(V\ot_{A}^LU, B)$. In the derived category of $A\ot B^{\rm op}$-modules we have 
\begin{multline*}
1_U\ot \beta=(1_U\ot \tilde\beta)(1_{U\ot_BV}\ot\alpha\ot 1_U)(1_{U\ot_BV\ot_AU}\ot \tilde\beta^{-1})\\
=(1_U\ot \tilde\beta)(\alpha\ot 1_{U\ot_BV\ot_AU})(1_{U\ot_BV\ot_AU}\ot \tilde\beta^{-1})
=(1_U\ot \tilde\beta)(\alpha\ot 1_U\ot \tilde\beta^{-1})=\alpha\ot 1_U.
\end{multline*}
Since $U\ot_B V\ot_A U$ is $A\ot B^{\rm op}$-projective, we have $\alpha\ot 1_U\sim 1_U\ot \beta$. Analogously, $1_V\ot\alpha\sim \beta\ot 1_V$.
\end{proof}

\begin{lemma}\label{techmain}
Suppose that $\alpha$ and $\beta$ satisfy the compatibility conditions of Lemma \ref{comp} and the maps
$$\ph_{\alpha\beta}:U\ot_B V\ot_A U\rightarrow U[1]\mbox{ and }\ph_{\beta\alpha}:V\ot_AU\ot_{B}V\rightarrow V[1]$$
are such that $\dd\ph_{\alpha\beta}=\alpha\ot 1_U-1_U\ot \beta$ and $\dd \ph_{\beta\alpha}=\beta\ot 1_V-1_V\ot\alpha$.
Then $$\beta(\ph_{\beta\alpha}\ot 1_U+1_V\ot \ph_{\alpha\beta}):V\ot_AU\ot_B V\ot_A U\rightarrow B[1]$$ is a null-homotopic chain map.
\end{lemma}
\begin{proof} Let us set $\psi=\beta(\ph_{\beta\alpha}\ot 1_U+1_V\ot\ph_{\alpha\beta})$. Since
\begin{multline*}
\psi d_{V\ot_AU\ot_B V\ot_A U}=\beta(\beta\ot 1_{V\ot_AU}-1_V\ot \alpha\ot 1_U+1_V\ot \alpha\ot 1_U-1_{V\ot_AU}\ot \beta)=0,
\end{multline*}
$\psi$ is a chain map. Note that $\psi(\beta\ot\beta)^{-1}\in\Hom_{\DDD^-B^e}(B,B[1])=0$. Consequently, $\psi$ equals $0$ in the derived category of $B^e$-modules. Since $V\ot_AU\ot_B V\ot_A U$ is projective, we have $\psi\sim 0$.
\end{proof}

Suppose that $A$ and $B$ are derived equivalent algebras, $U$ and $V$ are as above, and $\alpha$, $\beta$, $\ph_{\alpha\beta}$, and $\ph_{\beta\alpha}$ are as in Lemma \ref{techmain}.
If $(P,\mu_P)$ is a projective bimodule resolution of $A$, then it is easy to see that $\big(V\ot_A P\ot_A U,\beta(1_V\ot_A\mu_P\ot_A1_U)\big)=(\tilde P,\mu_{\tilde P})$ is a projective bimodule resolution of $B$.
For $f:P\rightarrow A[-n]$, we will denote by $\tilde f$ the map $1_V\ot_A f\ot_A 1_U:\tilde P\rightarrow B[-n]$.
There is an isomorphism $\chi:\HH^*(A)\rightarrow\HH^*(B)$ that sends the element corresponding to $f:P\rightarrow A[-n]$ to the element corresponding to $\chi(f)=\beta\tilde f:\tilde P\rightarrow B[-n]$. Note that $\mu_{\tilde P}=\chi(\mu_P)$.
Let now $\Delta_P:P\rightarrow P\ot_AP$ be a diagonal $2$-approximation for $(P,\mu_P)$. Since the map $1_P\ot_A\alpha\ot_A1_P:P\ot_A U\ot_BV\ot_A P\rightarrow P\ot_AP$ is a quasi-isomorphism and all the complexes under consideration are projective, there exists a chain map $\gamma:P\ot_AP\rightarrow P\ot_A U\ot_BV\ot_A P$ such that $\gamma(1_P\ot_A\alpha\ot_A1_P)\sim 1_{P\ot_A U\ot_BV\ot_A P}$ and $(1_P\ot_A\alpha\ot_A1_P)\gamma\sim 1_{P\ot_A P}$. Then it is easy to check that the map $\Delta_{\tilde P}=1_V\ot_A\gamma\Delta\ot_A1_U$ is a diagonal approximation for $(\tilde P,\mu_{\tilde P})$. Note also that $(1_P\ot_A\alpha\ot_A1_P)\gamma\Delta_P$ is a diagonal approximation for $(P,\mu_P)$. We have the following lemma.

\begin{lemma}\label{hlift} Let $f:P\rightarrow A[-n]$ be a map of $A^e$-complexes and $\phi_f:P\rightarrow P[1-n]$ be a homotopy lifting for $\big(f,(1_P\ot_A\alpha\ot_A1_P)\gamma\Delta_P\big)$. Then
$$
\psi_f=1_V\ot_A\phi_f\ot_A1_U+(-1)^n\big(\ph_{\beta\alpha}(\tilde f\ot_B1_V)\ot_A1_{P\ot_AU}+1_{V\ot_AP}\ot_A\ph_{\alpha\beta}(1_{U}\ot_B\tilde f)\big)\Delta_{\tilde P}
$$
is a homotopy lifting for $\big(\chi(f),\Delta_{\tilde P}\big)$.
\end{lemma}
\begin{proof} Direct calculations show that
\begin{multline*}
\dd\psi_f=(1_V\ot_Af\ot_A\alpha\ot_A1_{P\ot_AU}-1_{V\ot_AP}\ot_A\alpha\ot_Af\ot_A1_U)\Delta_{\tilde P}\\
+\big(\beta\tilde f\ot_B1_{\tilde P}-1_V\ot_Af\ot_A\alpha\ot_A1_{P\ot_AU}+1_{V\ot_AP}\ot_A\alpha\ot_Af\ot_A1_U-1_{\tilde P}\ot_B\beta\tilde f\big)\Delta_{\tilde P}\\
=\big(\chi(f)\ot_B1_{\tilde P}-1_{\tilde P}\ot_B\chi(f)\big)\Delta_{\tilde P}.
\end{multline*}
In particular, $\dd\psi_{\mu_P}=(\mu_{\tilde P}\ot_B1_{\tilde P}-1_{\tilde P}\ot_B\mu_{\tilde P})\Delta_{\tilde P}$. By the definition of the homotopy lifting, we have $\mu_P\phi_f+f\phi_{\mu_P}\sim 0$, and hence
\begin{multline*}
\mu_{\tilde P}\psi_f+\chi(f)\psi_{\mu_{P}}\sim (-1)^n\beta \big(\ph_{\beta\alpha}(\tilde f\ot_B1_V)\ot_A\mu_P\ot_A1_U+1_V\ot_A\mu_P\ot_A\ph_{\alpha\beta}(1_{U}\ot_B\tilde f)\big)\Delta_{\tilde P}\\
+\beta(1_V\ot_Af\ot_A1_U)\big(\ph_{\beta\alpha}(\tilde\mu_P\ot_B1_{V})\ot_A1_{P\ot_AU}+1_{V\ot_AP}\ot_A\ph_{\alpha\beta}(1_{U}\ot_B\tilde\mu_P)\big)\Delta_{\tilde P}\\
=(-1)^n\beta(\ph_{\beta\alpha}\ot_A1_U+1_V\ot_A\ph_{\alpha\beta})(\tilde f\ot_B\tilde\mu_P+\tilde\mu_P\ot_B\tilde f)\Delta_{\tilde P}\sim 0
\end{multline*}
by Lemma \ref{techmain}. Thus, $\psi_f$ is a homotopy lifting for $\big(\chi(f),\Delta_{\tilde P}\big)$.
\end{proof}

Now we are ready to prove the following theorem.

\begin{theorem}\label{DI} Suppose that $A$ and $B$ are $\kk$-algebras. If $A$ is derived equivalent to $B$, then $\HH^*(A)\cong\HH^*(B)$ as Gerstenhaber algebras.
\end{theorem}

\begin{proof} It is well known that the isomorphism $\chi$ defined above preserves the cup product. In fact it coincides with the isomorphism from \cite{Rickard}. Thus, it remains to prove that it preserves the Gerstenhaber bracket.

 By Lemma \ref{hlift} and Theorem \ref{main}, it is enough to show that
$$
(-1)^m\chi(f)\psi_g+(-1)^{m(n-1)}\chi(g)\psi_f\sim \chi\big((-1)^mf\phi_g+(-1)^{m(n-1)}g\phi_f\big)
$$
for any two maps $f:P\rightarrow A[-n]$ and $g:P\rightarrow A[-m]$ of $A^e$-complexes. We have
\begin{multline*}
(-1)^m\chi(f)\psi_g+(-1)^{m(n-1)}\chi(g)\psi_f-\chi\big((-1)^mf\phi_g+(-1)^{m(n-1)}g\phi_f\big)\\
=\beta\Big(1_V\ot_A\big((-1)^mf\phi_g+(-1)^{m(n-1)}g\phi_f\big)\ot_A1_U\Big)-\chi\big((-1)^mf\phi_g+(-1)^{m(n-1)}g\phi_f\big)\\
+\beta(1_V\ot_A f\ot_A1_U)\big(\ph_{\beta\alpha}(\tilde g\ot_B1_V)\ot_A1_{P\ot_AU}+1_{V\ot_AP}\ot_A\ph_{\alpha\beta}(1_{U}\ot_B\tilde g)\big)\Delta_{\tilde P}\\
-(-1)^{(m-1)(n-1)}\beta(1_V\ot_A g\ot_A1_U)\big(\ph_{\beta\alpha}(\tilde f\ot_B1_V)\ot_A1_{P\ot_AU}+1_{V\ot_AP}\ot_A\ph_{\alpha\beta}(1_{U}\ot_B\tilde f)\big)\Delta_{\tilde P}\\
=(-1)^{(m-1)n}\beta(\ph_{\beta\alpha}\ot_A1_U)(\tilde g\ot_B\tilde f)\Delta_{\tilde P}+(-1)^n\beta(1_V\ot_A\ph_{\alpha\beta})(\tilde f\ot_B\tilde g)\Delta_{\tilde P}\\
+(-1)^{n}\beta(\ph_{\beta\alpha}\ot_A1_U)(\tilde f\ot_B\tilde g)\Delta_{\tilde P}+(-1)^{(m-1)n}\beta(\ph_{\beta\alpha}\ot_A1_U)(\tilde f\ot_B\tilde g)\Delta_{\tilde P}\\
=(-1)^n\beta(\ph_{\beta\alpha}\ot_A1_U+1_V\ot_A\ph_{\alpha\beta})\big(\tilde f\ot_B\tilde g+(-1)^{mn}\tilde g\ot_B\tilde f\big)\Delta_{\tilde P}\sim 0
\end{multline*}
by Lemma \ref{techmain}. Thus, the theorem is proved.
\end{proof}

\section{A formula via contracting homotopy}\label{MAT}

In this section we present a formula that expresses the Lie bracket on the Hochschild cohomology in terms of an arbitrary resolution and a left contracting homotopy for it. Note that contracting homotopies can be used to construct the comparison maps between resolutions and this method was applied to compute the bracket, for example, in \cite{IIVZ}.

  \begin{definition} Let $(P,\mu_P)$ be a projective $A^e$-resolution of $A$. Let $t_P:P\rightarrow P$ and $\eta_P:A\rightarrow P$ be homomorphisms of  left modules such that $t_P(P_i)\subset P_{i+1}$ and $\eta_P(A)\subset P_0$.
  The pair $(t_P,\eta_P)$ is called a {\it left contracting homotopy} for $(P,\mu_P)$ if $d_Pt_P+t_Pd_P+\eta_P\mu_P=1_P$ and $t_P(t_P+\eta_P)=0$.
\end{definition}

Since $A$ is projective as a left $A$-module, any $A^e$-projective resolution of $A$ splits as a complex of left $A$-modules. Hence, a left contracting homotopy exists for any $A^e$-projective resolution of $A$ (see \cite[Lemma 2.3]{IIVZ} and the remark after it for details).

Let us fix an $A^e$-projective resolution $(P,\mu_P)$ of $A$ and a left contracting homotopy $(t_P,\eta_P)$ for it.

For any $n\ge 0$, the map $\pi_n:A\ot P_n\rightarrow P_n$ defined by the equality $\pi_n(a\otimes x)=ax$ for $a\in A$, $x\in P_n$ is an epimorphism of $A$-bimodules. Since $P_n$ is projective, there is $\iota_n\in\Hom_{A^e}(P_n, A\ot P_n)$ such that $\pi_n\iota_n=1_{P_n}$. Let us fix such $\iota_n$ for each $n\ge 0$. Then $\pi_n$ and $\iota_n$ ($n\ge 0$) determine homomorphisms of graded $A$-bimodules $\pi:A\ot P\rightarrow P$ and $\iota:P\rightarrow A\ot P$.

Let us define
$$
\begin{aligned}
t_{L}&:=(1_P\otimes\pi)(t_P\otimes 1_{P})(1_P\otimes\iota):P\otimes_A P\rightarrow (P\otimes_A P)[1],\\
 \eta_{L}&:=(1_P\otimes\pi)(\eta_P\otimes 1_{P})\iota:P\rightarrow P\otimes_A P,\\
 d_L&:=d_P\otimes 1_P, d_R:=1_P\otimes d_P:P\otimes_A P\rightarrow (P\otimes_A P)[-1],\\
 \mu_L&:=\mu_P\otimes 1_P,\mu_R:=1_P\otimes \mu_P:P\ot_AP\rightarrow P.
\end{aligned}
$$
Note that all the defined maps are homomorphisms of $A$-bimodules. Note also that we omit isomorphisms $\alpha_{P,P}^{ 1}$ and $\beta_{P,P}^{\pm 1}$ in our definitions according to our agreement. It is easy to see that the map $t_Ld_R:P\otimes_A P\rightarrow P\otimes_A P$ is locally nilpotent in the sense that for any $x\in P\otimes_A P$ there is an integer $l$ such that $(t_L d_R)^l(x)=0$. This follows from the fact that $t_L d_R(P\otimes_A P_j)\subset P\otimes_A P_{j-1}$ if $j>0$ and $t_L d_R(P\otimes_A P_0)=0$. Hence, the map $1_{P\otimes_A P}+ t_L d_R$ is invertible.

Let now $f:P\rightarrow A[-n]$ and $g:P\rightarrow A[-m]$ be maps of complexes. Let us define
$$
f\circ g=-f\mu_RS t_L(1_P\otimes g\otimes 1_P)(1_P\otimes S\eta_L)S\eta_L,
$$
where $S=(1_{P\otimes_A P}+ t_L d_R)^{-1}$.

\begin{theorem}\label{split_form} In the notation above the operation defined by the equality $[f,g]=f\circ g-(-1)^{(n-1)(m-1)}g\circ f$ induces the usual Lie bracket on the Hochschild cohomology.
\end{theorem}

We divide the proof into several lemmas. First of all, note that
\begin{equation}\label{eqs}
 d_L t_L+ t_L d_L+\eta_L\mu_L=1_{P\otimes_A P},\mu_L\eta_L=1_P, ( d_R)^2=( d_L)^2=0\mbox{ and } d_L d_R+ d_R d_L=0.
\end{equation}

\begin{lemma}\label{comm}
$( d_L+ d_R)S=S( d_L+\eta_L\mu_L d_R)$.
\end{lemma}
\begin{proof}
Let us multiply the desired equality by $1_{P\otimes_AP}+ t_L d_R$ on the left and on the right at the same time. We obtain that we have to prove that
$$
 d_L+ d_R+ t_L d_R d_L+ t_L( d_R)^2= d_L+\eta_L\mu_L d_R+ d_L t_L d_R+\eta_L\mu_L d_R t_L d_R.
$$
Using \eqref{eqs} one can see that it is enough to show that $\eta_L\mu_L d_R t_L d_R=0$. But the last equality follows from the fact that the image of $ d_R t_L d_R$ lies in $\oplus_{n>0}P_n\otimes_A P\subset\Ker\mu_L$.
\end{proof}

\begin{lemma}\label{appr}
$S\eta_L$ is a diagonal $2$-approximation of $P$.
\end{lemma}
\begin{proof} By Lemma \ref{comm} we have $$\dd(S\eta_L)=( d_L+ d_R)S\eta_L-S\eta_Ld_P=S( d_L+\eta_L\mu_L d_R)\eta_L-S\eta_Ld_P.$$
Since $\Im\eta_L\subset\Ker d_L$, it is enough to prove that $\eta_L\mu_L d_R\eta_L=\eta_Ld_P$. It is easy to see that $\mu_L d_R=d_P\mu_L$. Hence,
$\eta_L\mu_L d_R\eta_L=\eta_Ld_P\mu_L\eta_L=\eta_Ld_P$ by \eqref{eqs}.
\end{proof}

\begin{proof}[Proof of Theorem \ref{split_form}.] It follows from Lemma \ref{appr} that 
$\Delta_P=(1_P\otimes \mu_RS \eta_L)S\eta_L$ is a diagonal $2$-approximation of $P$.

It is enough to show that
$
\phi_g=(-1)^{m-1}\mu_RS t_L(1_P\otimes g\otimes 1_P)(1_P\otimes S\eta_L)S\eta_L
$
is a homotopy lifting for $(g,\Delta_P)$. Using Lemma \ref{comm}, we get $\mu_R\dd(S)t_L=\mu_R(1_P-\eta_L\mu_L)d_Rt_L=\mu_Rd_Rt_L-\mu_R\eta_Ld_P\mu_Lt_L=0$.
Since $Sd_R=d_R$ and $\mu_LS\eta_L=1_P$, we get now
\begin{multline*}
\dd\phi_g=-\mu_RS\big((d_L+d_R)t_L+t_L(d_L+d_R)\big)(1_P\otimes g\otimes 1_P)(1_P\otimes S\eta_L)S\eta_L\\
=\mu_RS(\eta_L\mu_L-1_{P\ot_AP}-t_Ld_R)(1_P\otimes g\otimes 1_P)(1_P\otimes S\eta_L)S\eta_L-\mu_Rd_Rt_L(1_P\otimes g\otimes 1_P)(1_P\otimes S\eta_L)S\eta_L\\
=(\mu_P\otimes g\otimes \mu_RS\eta_L)(1_P\otimes S\eta_L)S\eta_L-(1_P\otimes g\otimes \mu_P)(1_P\otimes S\eta_L)S\eta_L\\
=(g\otimes 1_P)(1_P\otimes \mu_RS\eta_L)S\eta_L\mu_LS\eta_L-(1_P\otimes g)(1_P\otimes \mu_RS\eta_L)S\eta_L
=(g\otimes 1_P-1_P\otimes g)\Delta_P.
\end{multline*}
Note also that $\mu_L\Delta_P=\mu_RS \eta_L=\mu_R\Delta_P$ and $\mu_P\phi_g=0$. Hence, $\phi_g$ is a homotopy lifting for $(g,\Delta_P)$ and the theorem is proved.
\end{proof}

\section{Formulas for the Connes' differential}\label{WTP}

In this section we discuss some formulas for the Connes' differential. These formulas are based on the formula from \cite{Kaled}. In the case of a symmetric algebra a formula for the Connes' differential gives a formula for a BV differential. Thus, we obtain in this section an alternative way for computing the Lie bracket on the Hochschild cohomology of a symmetric algebra.

Let $Tr$ denote the functor $A\ot_{A^e}-$ from the category of $A$-bimodules to the category of $\kk$-linear spaces. If $M$ and $N$ are $A$-bimodules, then there is an isomorphism $\si_{M,N}:Tr(M\ot_AN)\rightarrow Tr(N\ot_AM)$ defined by the equality $\si_{M,N}(1\ot x\ot y)=1\ot y\ot x$ for $x\in M$ and $y\in N$. Moreover, for $f\in\Hom_{A^e}(M_1,M_2)$ and $g\in\Hom_{A^e}(N_1,N_2)$ one has $\si_{M_2,N_2}Tr(f\ot g)=Tr(g\ot f)\si_{M_1,N_1}$. It is easy to see also that $Tr$ induces a functor from the category of $A^e$-complexes to the category of $\kk$-complexes. In this case $\si_{P,Q}$ is defined by the equality $\si_{P,Q}(1\ot x\ot y)=(-1)^{ij}\ot y\ot x$ for $x\in P_i$ and $y\in Q_j$ and satisfies the property $\si_{P_2,Q_2}Tr(f\ot g)=Tr(g\ot f)\si_{P_1,Q_1}$ for $f:P_1\rightarrow P_2$ and $g:Q_1\rightarrow Q_2$.

The Hochschild homology $\HH_*(A)$ of the algebra $A$ is simply the homology of the complex $Tr(\Bar(A))$. As in the case of cohomology, any comparison morphism $\Phi_P^Q:P\rightarrow Q$ between resolutions $(P,\mu_P)$ and $(Q,\mu_Q)$ of the algebra $A$ induces an isomorphism $Tr(\Phi_P^Q):\Ho_*Tr(P)\rightarrow \Ho_*Tr(Q)$. Thus, the Hochschild homology of $A$ is isomorphic to the homology of $Tr(P)$ for any projective bimodule resolution $(P,\mu_P)$ of $A$.

Note that $Tr(\Bar_n(A))\cong A^{\ot (n+1)}$. Connes' differential $\BB:\HH_n(A)\rightarrow\HH_{n+1}(A)$ is the map induced by the map from $Tr(\Bar_n(A))$ to $Tr(\Bar_{n+1}(A))$ that sends $a_0\ot a_1\ot\dots\ot a_n\in A^{\ot(n+1)}$ to
$$
\sum_{i=0}^n (-1)^{in} 1\otimes a_i\otimes\dots \otimes
a_n\otimes a_0\otimes \cdots \otimes a_{i-1}+ \sum_{i=0}^n
(-1)^{in}a_i\otimes 1 \otimes  a_{i+1}\otimes \dots \otimes
a_n\otimes a_0\otimes \dots \otimes a_{i-1}.
$$
In fact, it follows from some standard arguments that the homological class of the second summand is zero. The following result is essentially stated in \cite{Kaled} (see equation (4.8) of the cited paper and the explanation before and after it).

\begin{proposition}[D. Kaledin]\label{Kaledin}
Let $(P,\mu_P)$ be a projective bimodule resolution of $A$, $\Delta_P$ be a diagonal $2$-approximation for $P$, and $\phi_P:P\ot_AP\rightarrow P[1]$ be such that $\mu_P\otimes 1_P-1_P\otimes\mu_P=\dd\phi$. Then the map $$Tr(\phi_P)(1_{P\ot_AP}+\si_{P,P})Tr(\Delta_P):Tr(P)\rightarrow Tr(P[1])$$ induces the Connes's differential on the Hochschild homology.
\end{proposition}

This result can be written in a slightly different form.

\begin{corollary}\label{CorKaled}
Let $(P,\mu_P)$, $\Delta_P$, and $\phi_P$ be as in Proposition \ref{Kaledin}, and $\ee:P\rightarrow P[1]$ be such that $(\mu_P\otimes 1_P-1_P\otimes\mu_P)\Delta_P=\dd\ee$. Then the map $$Tr(\phi_P)\si_{P,P}Tr(\Delta_P)+Tr(\ee):Tr(P)\rightarrow Tr(P[1])$$ induces the Connes's differential on the  Hochschild homology.
\end{corollary}
\begin{proof} Since $\dd(\phi_P\Delta_P)=(\mu_P\otimes 1_P-1_P\otimes\mu_P)\Delta_P$, it is enough to note that the map $H_*(Tr(\phi_P)\si_{P,P}Tr(\Delta_P)+Tr(\ee)):\HH_*(A)\rightarrow\HH_*(A)$ does not depend on the choice of $\ee$.
\end{proof}

Now it is not difficult to express the Connes' differential in terms of a contracting homotopy.

\begin{corollary}\label{DifMAT} Let $S$, $t_L$ and $\eta_L$ be as in the previous section. Then the map $-Tr(\mu_RSt_L)\si_{P,P}Tr\big((1_P\ot (\mu_RS\eta_L)^2)S\eta_L\big)$ induces the Connes' differential on the Hochschild homology.
\end{corollary}
\begin{proof} It follows from Lemma \ref{appr} that $\mu_RS\eta_L:P\rightarrow P$ is a comparison morphism, i.e. there is some $u:P\rightarrow P[1]$ such that $1-\mu_RS\eta_L=\dd u$. It is not hard to show using Lemma \ref{comm} (see also the proof of Theorem \ref{split_form}) that
$\dd\phi_P=\mu_L-\mu_R$ for $\phi_P=u(\mu_L-\mu_R)-\mu_RSt_L(\mu_RS\eta_L\ot 1_P)$. Let also note that $(\mu_L-\mu_R)\Delta_P=0$ for $\Delta_P=(1_P\ot\mu_RS\eta_L)S\eta_L$. Then the Connes' differential is induced by the map
$$
Tr(\phi_P)\si_{P,P}Tr(\Delta_P)=-Tr(\mu_RSt_L)\si_{P,P}Tr\big((1_P\ot (\mu_RS\eta_L)^2)S\eta_L\big).
$$
\end{proof}

Now we explain how one can obtain a formula for a BV differential on the Hochschild cohomology of a symmetric algebra in terms of an arbitrary resolution.

First of all, let us recall that there are well known maps $\ii_f:\HH_*(A)\rightarrow \HH_*(A)$  for $f\in\HH^*(A)$, whose definition can be found, for example, in \cite{Menichi}. These maps satisfy the condition $\ii_f\ii_g=\ii_{f\smile g}$. We need also the fact that $\ii_f|_{\HH_n(A)}=0$ for $n<|f|$ and that
 $\ii_f|_{\HH_{|f|}(A)}$ is the map induced by $Tr(\tilde f):Tr(P_n)\rightarrow Tr(A)\cong\HH_0(A)$, where $\tilde f\in\Hom_{A^e}(P_n,A)$ represents $f$. After the correction of signs one obtains by \cite[Lemma 15]{Menichi} that
\begin{multline*}
\ii_{[f,g]}(x)=(-1)^{(|f|+1)|g|}(-(-1)^{|f|+|g|)}\BB\ii_{f\smile g}(x)+\ii_f \BB\ii_g(x)-(-1)^{|f||g|}\ii_g \BB\ii_f( x)-\ii_{f\smile g}\BB(x))
\end{multline*}
for all $f,g\in\HH^*(A)$, $x\in\HH_*(A)$. Considering $x\in\HH_{|f|+|g|-1}(A)$, we get
\begin{equation}\label{Tr}
Tr([f,g])=-(-1)^{(|f|+1)|g|}(Tr(f\smile g)\BB-Tr(f)\BB\ii_g-(-1)^{|f||g|}Tr(g) \BB\ii_f).
\end{equation}

\begin{definition} 
A {\it Batalin--Vilkovisky  algebra} ({\it BV algebra} for short) is a
Gerstenhaber algebra $(R^\bullet,\, \smile,\, [\ \,,\,\ ])$
with an operator $\DD\colon R^\bullet \rightarrow
R^{\bullet-1}$ of degree $-1$ such that $\DD\circ  \DD=0$
and
$$
[a,\, b]=-(-1)^{(|a|+1)|b|}(\DD(a\smile b)- \DD(a)\smile
b-(-1)^{|a|}a\smile \DD(b))
$$
for   homogeneous  elements $a, b\in R^\bullet$.
\end{definition}

\begin{definition}
The finite dimensional algebra $A$ is called {\it symmetric} if $A\cong \Hom_{\kk}(A,k)$ as an $A$-bimodule.
\end{definition}

Let $A$ be symmetric. Let $\theta:A\rightarrow\kk$ be an image of $1$ under some bimodule isomorphism from $A$ to $\Hom_{\kk}(A,\kk)$. Then it is easy to see that $\theta$ induces a map from $Tr(A)$ to $\kk$. We  denote this map by $\theta$ too. Note also that if $f\in\Hom_{A^e}(M,A)$, then $\theta Tr(f)=0$ if and only if $f=0$.

Let $\BB_P:Tr(P)\rightarrow Tr(P[1])$ be a map inducing the Connes' differential on the Hochschild homology. Then we can  define $\DD_P(f):P\rightarrow A[1-n]$ for $f:P\rightarrow A[-n]$ as the unique map such that $\theta Tr(\DD_P(f))=\theta Tr(f)\BB_P$.

\begin{proposition}[\cite{Tradler}]\label{BVT} $\DD_P$ induces  a BV differential on the Hochschild cohomology.
\end{proposition}

Proposition \ref{BVT} is the remark after \cite[Theorem 1]{Tradler}. To see that it is valid, one can apply $\theta$ to the equality \eqref{Tr} with $\BB=\BB_P$ and get
\begin{multline*}
\theta Tr([f,g])=-(-1)^{(|f|+1)|g|}(\theta Tr \DD_P(f\smile g)-\theta Tr\DD_P(f)\ii_g-(-1)^{|f||g|}\theta Tr\DD_P(g) \ii_f)\\
=-(-1)^{(|f|+1)|g|}\theta Tr( \DD_P(f\smile g)-\DD_P(f)\smile g-(-1)^{|f|} f\smile \DD_P(g)).
\end{multline*}
Note also that if one knows the BV differential and the cup product, then it is easy to compute the Gerstenhaber bracket.

\section{Example of an application}\label{EX}

In this section we apply the results of the previous sections to describe the BV structure on the Hochschild cohomology of the family of algebras  considered in \cite{Gen1} and \cite{Gen2}.
During this section we fix some integer $k>1$ and set $A=\kk\la x_0,x_1\ra/\la x_0^2,x_1^2,(x_0x_1)^k-(x_1x_0)^k\ra$. The index $\alpha$ in the notation $x_{\alpha}$ is always specified modulo $2$. If $a$ is an element of $\kk\la x_0,x_1\ra$, then we denote by $a$ its class in $A$ too.

Let $G$ be a subset of $\kk\la x_0,x_1\ra$ formed by the elements $(x_0x_1)^{i+1}$, $x_1(x_0x_1)^i$, $(x_1x_0)^i$, and $x_0(x_1x_0)^i$ for $0\le i\le k-1$.
Note that the classes of the elements from $G$ form a basis of $A$. Let $G$ denote this basis too. Let $l_v$ denote the length of $v\in G$.
Note that the algebra $A$ is symmetric with $\theta$ defined by the equalities $\theta\big((x_0x_1)^k\big)=1$ and $\theta(v)=0$ for $v\in G\setminus\{(x_0x_1)^k\}$.
For $v\in G$, we introduce $v^*\in G$ as the unique element such that $\theta(vv^*)=1$. Note that $\theta(vw)=0$ for $w\in G\setminus\{v^*\}$. For $a=\sum\limits_{v\in G}a_vv\in A$, where $a_v\in \kk$ for $v\in G$, we define $a^*:=\sum\limits_{v\in G}a_vv^*\in A$.
It is clear that $(a^*)^*=a$ for any $a\in A$. If $v,w\in G$, then $\frac{v}{w}$ denotes $(v^*w)^*$. If there is such $u\in G$ that $wu=v$, then this $u$ is unique and $\frac{v}{w}=u$. If there is no such $u$, then $\frac{v}{w}=0$. Note that $\frac{\frac{v}{x_{\alpha}}x_{\beta}}{x_{\alpha}}$ is equal to $\frac{v}{x_{\alpha}}$ if $\alpha=\beta$ and $v\in\{x_{\alpha},1^*\}$, and is equal to $0$ in all remaining cases.
 For $a=\sum\limits_{v\in G}a_vv\in A$ and $b=\sum\limits_{v\in G}b_vv\in A$, where $a_v,b_v\in \kk$ for $v\in G$, we define $\frac{a}{b}:=\sum\limits_{v,w\in G}a_vb_w\frac{v}{w}\in A$.

In this section we will use the bimodule resolution of $A$ described in \cite{Gen1}. Here we present it in a little another form, but one can easily check that it is the same resolution.
Let us introduce the algebra $B=\kk[x_0,x_1,z]/\la x_0x_1\ra$. We introduce the grading on $B$ by the equalities $|x_0|=|x_1|=1$ and $|z|=2$. Let us define the $A^e$-complex $P$.
We set $P=A\ot B\ot A$ as an $A$-bimodule. The grading on $P$ comes from the grading on $B$ and the trivial grading on $A$. Let $\und{a}$ ($a\in B$) denote $1\ot a\ot 1$. For convenience we set $\und{a}=0$ if $a=x_{\alpha}^iz^j$, where ${\alpha}\in\{0,1\}$ and $i$ or $j$ is less than $0$.
We define the differential $d_P$ by the equality
$$
d_P(\und{x_{\alpha}^iz^j})=\begin{cases}
0,&\mbox{ if $i=j=0$},\\
x_{\alpha} \und{x_{\alpha}^{i-1}}+(-1)^i\und{x_{\alpha}^{i-1}} x_{\alpha},&\mbox{ if $j=0$, $i>0$},\\
\sum\limits_{v\in G,{\beta}\in\{0,1\}}(-1)^{jl_v+{\beta}}v^*\und{x_{\beta}z^{j-1}}\frac{v}{x_{\beta}},&\mbox{ if $i=0$, $j>0$},\\
x_{\alpha}\und{x_{\alpha}^{i-1}z^j}+(-1)^{i+j}\und{x_{\alpha}^{i-1}z^j} x_{\alpha}\\
+(-1)^{i+{\alpha}}((-1)^jx_{\alpha}^*\und{x_{\alpha}^{i+1}z^{j-1}}+\und{x_{\alpha}^{i+1}z^{j-1}} x_{\alpha}^*),&\mbox{ if $i,j>0$}.
\end{cases}
$$
for ${\alpha}\in\{0,1\}$, $i,j\ge 0$. We define $\mu_P:P_0\rightarrow A$ by the equality $\mu_P(\und{1})=1$. Then one can check that $(P,\mu_P)$ is an $A^e$-projective resolution of $A$ isomorphic to the resolution from \cite{Gen1}. Let us define the left contracting homotopy $(t_P,\eta_P)$ for $(P,\mu_P)$. We define $\eta_P$ by the equality $\eta_P(1)=\und{1}$. Now, for $v\in G$, ${\alpha}\in\{0,1\}$ and $i,j\ge 0$, we define
$$
t_P(\und{x_{\alpha}^iz^j}v)=\begin{cases}
\sum\limits_{w\in G,{\beta}\in\{0,1\}}(-1)^{j(l_w+l_v+1)+1}\frac{w^*}{v^*}\und{x_{\beta}z^j}\frac{w}{x_{\beta}},&\mbox{ if $i=0$, $v\not=1^*$},\\
\sum\limits_{w\in G}(-1)^{j(l_w+1)+1}w^*\und{x_{l_w}z^j}\frac{w}{x_{l_w}},&\mbox{ if $i=0$ and $v=1^*$},\\
(-1)^{i+j+1}\und{x_{\alpha}^{i+1}z^j}\frac{v}{x_{\alpha}}+(-1)^{jl_v+j+l_v}\frac{v}{x_{1}^*}\und{z^{j+1}},&\mbox{ if $i=1$ and ${\alpha}=1$},\\
(-1)^{i+j+1}\und{x_{\alpha}^{i+1}z^j}\frac{v}{x_{\alpha}},&\mbox{ otherwise}.
\end{cases}
$$

In this section we will use the notation of Section \ref{MAT}. Our aim is to describe the BV structure on the Hochschild cohomology of $A$. As it was explained in the previous section, it is enough to describe the Connes' differential.
By Corollary \ref{DifMAT} we have to describe the map $-Tr(St_L)\si_{P,P}Tr\big((1_P\ot (\mu_RS\eta_L)^2)S\eta_L\big)$. Let us start with the map $S\eta_L:P\rightarrow P\ot_A P$.

Firstly, let introduce the following notation:
$$
\begin{aligned}
A_{t,j}&=\sum\limits_{\tiny\begin{array}{c}v,w\in G\\ \alpha,\beta\in\{0,1\}\end{array}}(-1)^{jl_v+t(l_w+l_v+1)+{\beta}}\frac{w^*}{v}\und{x_{\alpha}z^t}\frac{w}{x_{\alpha}}\ot\und{x_{\beta}z^{j-1}}\frac{v}{x_{\beta}},\\
B_{t,j}&=\sum\limits_{v\in G, {\beta}\in\{0,1\}}(-1)^{(j+t)(l_v+1)}\frac{x_{{\beta}+1}^*}{v}\und{x_{{\beta}+1}^2z^t}(x_{\beta}x_{{\beta}+1})^{k-1}\ot\und{x_{\beta}^2z^{j-2}}\frac{v}{x_{\beta}},\\
C_{t,i,j,\alpha}&=(-1)^{i+j+{\alpha}}\sum\limits_{w\in G,\beta\in\{0,1\}}(-1)^{tl_w}\frac{w^*}{x_{\alpha}}\und{x_{\beta}z^t}\frac{w}{x_{\beta}}\ot \und{x_{\alpha}^{i+1}z^{j-1}},\\
D_{t,i,j,\alpha}&=(-1)^{(i+1)t}\sum\limits_{v\in G,{\beta}\in\{0,1\}}(-1)^{jl_v+{\beta}}\und{x_{\alpha}^{i+1}z^t}\frac{v^*}{x_{\alpha}}\ot\und{x_{\beta}z^{j-1}}\frac{v}{x_{\beta}},\\
E_{t,i,j,\alpha}&=(x_{{\alpha}+1}x_{\alpha})^{k-1}\und{x_{{\alpha}+1}^2z^t}(x_{\alpha}x_{{\alpha}+1})^{k-1}\ot \und{x_{\alpha}^{i+2}z^{j-2}}\\
&\phantom{=(x_{{\alpha}+1}x_{\alpha})^{k-1}\und{x_{{\alpha}+1}^2z^t}(x_{\alpha}x_{{\alpha}+1})^{k-1}}+(-1)^{it}\und{x_{{\alpha}}^{i+2}z^t}(x_{{\alpha}+1}x_{{\alpha}})^{k-1}\ot\und{x_{{\alpha}+1}^2z^{j-2}}(x_{{\alpha}}x_{{\alpha}+1})^{k-1}.
\end{aligned}
$$

\begin{lemma}\label{zt} If $q,j\ge 0$ are some integers, then
$$
S(\und{z^q}\ot \und{z^j})=\sum\limits_{t=0}^j(-1)^{(j+q+1)t}\big(\und{z^{q+t}}\ot\und{z^{j-t}}+A_{q+t,j-t}+B_{q+t,j-t}\big).
$$
In particular, $S\eta_L(\und{z^j})=\sum\limits_{t=0}^j(-1)^{(j+1)t}\big(\und{z^t}\ot\und{z^{j-t}}+A_{t,j-t}+B_{t,j-t}\big).$
\end{lemma}
\begin{proof} We have to show that
\begin{equation}\label{toshow}
\und{z^q}\ot\und{z^j}=\sum\limits_{t=0}^j(-1)^{(j+q+1)t}(1_{P\ot_AP}+t_Ld_R)(\und{z^{q+t}}\ot\und{z^{j-t}}+A_{q+t,j-t}+B_{q+t,j-t}).
\end{equation}
Direct calculations show that $t_Ld_R(\und{z^{q+t}}\ot\und{z^{j-t}})+A_{q+t,j-t}=0$ and $t_Ld_RB_{q+t,j-t}=0$ for $0\le t\le j$.
One can show that if $t_P\big(\und{x_{\alpha}z^{q+t}}\frac{w}{x_{\alpha}}x_{\beta}\big)\not=0$, then either $\frac{w^*}{x_{\beta}}=0$ or $w=x_0^*$, $\alpha=1$ and $\beta=0$. In the first case $\frac{w^*}{v}=0$ or $\frac{v}{x_{\beta}}=0$ for any $v\in G$. In the second case we have $t_P\big(\und{x_{1}z^{q+t}}\frac{x_0^*}{x_{1}}x_{0}\big)=-\und{z^{q+t+1}}$, and $\frac{x_0}{v}\not=0$ and $\frac{v}{x_{0}}\not=0$ simultaneously only for $v=x_0$. 
Analogously, we have $t_P\big(\und{x_{\alpha}z^{q+t}}\frac{w}{x_{\alpha}}x_{\beta}^*\big)\not=0$ only if either $\frac{w^*}{x_{\beta}}=0$ or $w=x_{\alpha}$, $\alpha=\beta+1$. In the last case $t_P(\und{x_{\alpha}z^{q+t}}x_{\beta}^*)=(-1)^{q+t}\und{x_{{\beta}+1}^2z^{q+t}}(x_{\beta}x_{{\beta}+1})^{k-1}$. Thus,
$$
t_Ld_RA_{q+t,j-t}+B_{q+t,j-t}+(-1)^{j+q+1}\und{z^{q+t+1}}\ot\und{z^{j-t-1}}=0.
$$ Substituting the obtained values of $t_Ld_R(\und{z^{q+t}}\ot\und{z^{j-t}})$, $t_Ld_RA_{q+t,j-t}$, and $t_Ld_RB_{q+t,j-t}$ to \eqref{toshow} we obtain a true equality.
\end{proof}

\begin{lemma}\label{uizt} If ${\alpha}\in\{0,1\}$, and $q,j\ge0$ and $i>0$ are some integers, then
\begin{multline*}
S(\und{z^q}\ot \und{x_{\alpha}^iz^j})=\sum\limits_{t=0}^j(-1)^{(i+j+q+1)t}\Big(\sum\limits_{r=0}^i(-1)^{r(q+t)}\und{x_{\alpha}^rz^{q+t}}\ot\und{x_{\alpha}^{i-r}z^{j-t}}\\+C_{q+t,i,j-t,\alpha}+D_{q+t,i,j-t,\alpha}+E_{q+t,i,j-t,\alpha}\Big).
\end{multline*}
In particular,
$$
S\eta_L(\und{x_{\alpha}^iz^j})=\sum\limits_{t=0}^j(-1)^{(i+j+1)t}\Big(\sum\limits_{r=0}^i(-1)^{rt}\und{x_{\alpha}^rz^t}\ot\und{x_{\alpha}^{i-r}z^{j-t}}+C_{t,i,j-t,\alpha}+D_{t,i,j-t,\alpha}+E_{t,i,j-t,\alpha}\Big).
$$
\end{lemma}
\begin{proof} Firstly, note that $t_P(\und{x_{\alpha}^rz^{q+t}}x_{\alpha})=(-1)^{r+q+t+1}\und{x_{\alpha}^{r+1}z^{q+t}}$. Also $t_P(\und{x_1z^{q+t}}x_1^*)=-\und{z^{q+t+1}}$ and $t_P(\und{x_{\alpha}^rz^{q+t}}x_{\alpha}^*)=0$ for $r>1$ and for $r=1$, ${\alpha}=0$. Hence, we have
\begin{multline*}
(1_{P\ot_AP}+t_Ld_R)\left(\sum\limits_{r=0}^i(-1)^{r(q+t)}\und{x_{\alpha}^rz^{q+t}}\ot\und{x_{\alpha}^{i-r}z^{j-t}}\right)\\
=\und{z^{q+t}}\ot\und{x_{\alpha}^{i}z^{j-t}}-C_{q+t,i,j-t,\alpha}-D_{q+t,i,j-t,\alpha}+(-1)^{i+j+q}\alpha\und{z^{q+t+1}}\ot\und{x_{\alpha}^iz^{j-t-1}},
\end{multline*}
Now we have $t_P\big(\und{x_{\beta}z^{q+t}}\frac{w}{x_{\beta}}x_{\alpha}^*\big)\not=0$ only if either $w=x_{\beta}$ or $w=x_{\alpha+1}x_{\alpha}$, $\beta=\alpha+1$. In the last case and in the case $w=x_{\beta}$, $\beta=\alpha$ we have $\frac{w^*}{x_{\alpha}}=0$. Thus, the only nonzero case is  $t_P(\und{x_{\alpha+1}z^{q+t}}x_{\alpha}^*)=(-1)^{q+t}\und{x_{{\alpha}+1}^2z^{q+t}}(x_{{\alpha}}x_{{\alpha}+1})^{k-1}$. Further, we have $t_P\big(\und{x_{\beta}z^{q+t}}\frac{w}{x_{\beta}}x_{\alpha}\big)\not=0$ only if either $\frac{w^*}{x_{\alpha}}=0$ or $w=x_0^*$, $\alpha=0$ and $\beta=1$. In the last case we have  $t_P\big(\und{x_{1}z^{q+t}}x_1^*\big)=-\und{z^{q+t+1}}$. Thus,
\begin{multline*}
t_Ld_RC_{q+t,i,j-t,\alpha}\\=-(x_{{\alpha}+1}x_{\alpha})^{k-1}\und{x_{{\alpha}+1}^2z^{q+t}}(x_{{\alpha}}x_{{\alpha}+1})^{k-1}\ot \und{x_{\alpha}^{i+2}z^{j-t-2}}+(-1)^{i+j+q}(1-\alpha)\und{z^{q+t+1}}\ot\und{x_{\alpha}^iz^{j-t-1}}.
\end{multline*}
One can check that $t_P\big(\und{x_{\alpha}^{i+1}z^{q+t}}\frac{v^*}{x_{\alpha}}x_{\beta}\big)=0$ if $\frac{v}{x_{\beta}}\not=0$. If $\frac{v^*}{x_{\alpha}}x_{\beta}^*\not=0$, then either $v=x_{\alpha}^*$ or $v=(x_{\alpha}x_{{\alpha}+1})^*$, $\beta=\alpha+1$. Since in the second case $\frac{v}{x_{\beta}}=0$, we have
\begin{multline*}
t_Ld_RD_{q+t,i,j-t,\alpha}=(-1)^{(i+1)(q+t+1)}t_L\left(\sum\limits_{{\beta}\in\{0,1\}}\und{x_{\alpha}^{i+1}z^{q+t}}x_{\beta}^*\ot\und{x_{\beta}^2z^{j-t-2}}\frac{x_{\alpha}^*}{x_{\beta}}\right)\\
=(-1)^{i(q+t)+1}\und{x_{\alpha}^{i+2}z^{q+t}}(x_{{\alpha}+1}x_{\alpha})^{k-1}\ot\und{x_{{\alpha}+1}^2z^{j-t-2}}(x_{\alpha}x_{{\alpha}+1})^{k-1}.
\end{multline*}
Finally, note that $t_Ld_RE_{q+t,i,j-t,\alpha}=0.$
Taking in account all the proved equalities, we obtain
\begin{multline*}
\sum\limits_{t=0}^j(-1)^{(i+j+q+1)t}(1_{P\ot_AP}+t_Ld_R)\Big(\sum\limits_{r=0}^i(-1)^{rt}\und{x_{\alpha}^rz^{q+t}}\ot\und{x_{\alpha}^{i-r}z^{j-t}}+C_{q+t,i,j-t,\alpha}+D_{q+t,i,j-t,\alpha}\\
+E_{q+t,i,j-t,\alpha}\Big)
=\sum\limits_{t=0}^j(-1)^{(i+j+q+1)t}(\und{z^{q+t}}\ot\und{x_{\alpha}^{i}z^{j-t}}+(-1)^{i+j+q}\und{z^{q+t+1}}\ot\und{x_{\alpha}^iz^{j-t-1}})=\und{z^q}\ot\und{x_{\alpha}^{i}z^{j}}.
\end{multline*}
\end{proof}

From Lemmas \ref{zt} and \ref{uizt} we obtain the following statement.

\begin{lemma}\label{zz} $\mu_RS(\und{z^q}\ot\und{x_{\alpha}^rz^t})=(-1)^{q(r+t)}\und{x_{\alpha}^rz^{q+t}}$. In particular, $\mu_RS\eta_L=1_P$.
\end{lemma}

It remains to describe $Tr(\mu_RSt_L)\si_{P,P}$ on the image of $Tr(S\eta_L)$.

\begin{lemma}\label{zero}
Let $v\in G$, ${\alpha},\beta\in\{0,1\}$, $p$, $r$, $q$ and $t$ be some integers. Suppose that $p>0$ and  one of the conditions $p=1$, ${\alpha}=1$, and $v\in\{x_1^*,1^*\}$ is not satisfied. Then 
\begin{multline*}
\mu_RSt_L(\und{x_{\alpha}^pz^q}v\ot\und{x_{\beta}^rz^t})\\
=\begin{cases}
(-1)^{(r+1)q+p+1}\und{x_{\alpha}^{p+r+1}z^q}\frac{v}{x_{\alpha}},&\mbox{if $t=0$ and either $v\in\{x_{\alpha},1^*\}$, $\beta=\alpha$ or $r=0$},\\
0,&\mbox{otherwise}.
\end{cases}
\end{multline*}
\end{lemma}
\begin{proof} Assume firstly that $\beta\not=\alpha$ and $r>0$. Direct calculations show that $(t_Ld_R)^2t_L(\und{x_{\alpha}^pz^q}v\ot\und{x_{\beta}^rz^t})=0$. Then
$$
\mu_RSt_L(\und{x_{\alpha}^pz^q}v\ot\und{x_{\beta}^rz^t})=\mu_R(t_L-t_Ld_Rt_L+S(t_Ld_R)^2t_L)(\und{x_{\alpha}^pz^q}v\ot\und{x_{\beta}^rz^t})=0.
$$
Here we also use the fact that $t_Ld_Rt_L(\und{x_{\alpha}^pz^q}v\ot\und{x_{\beta}})=0$ for the case $r=1$, $t=0$. In the remaining part of the proof we assume that $\beta=\alpha$.

Let us consider the case where $r>0$. Direct calculations show that $t_Ld_Rt_L(\und{x_{\alpha}^pz^q}v\ot\und{x_{{\alpha}}^rz^t})=0$ if $v\in G\setminus\{x_{\alpha},1^*\}$. One can also check that
$t_Ld_Rt_L(\und{x_{\alpha}^pz^q}v\ot\und{x_{{\alpha}}^rz^t})=(-1)^qt_L(\und{x_{\alpha}^{p+1}z^q}v\ot\und{x_{{\alpha}}^{r-1}z^t})$ if $v\in \{x_{\alpha},1^*\}$. Then
\begin{multline*}
\mu_RSt_L(\und{x_{\alpha}^pz^q}v\ot\und{x_{{\alpha}}^rz^t})=\mu_R(t_L-St_Ld_Rt_L)(\und{x_{\alpha}^pz^q}v\ot\und{x_{{\alpha}}^rz^t})\\
=\begin{cases}
(-1)^{q+1}\mu_RSt_L(\und{x_{\alpha}^{p+1}z^q}v\ot\und{x_{{\alpha}}^{r-1}z^t}),&\mbox{if $v\in\{x_{\alpha},1^*\}$},\\
0,&\mbox{otherwise}.
\end{cases}
\end{multline*}
If $v\in G\setminus\{x_{\alpha},1^*\}$, then the required equality is proved. If $v\in\{x_{\alpha},1^*\}$, then we obtain 
$\mu_RSt_L(\und{x_{\alpha}^pz^q}v\ot\und{x_{{\alpha}}^rz^t})=(-1)^{r(q+1)}\mu_RSt_L(\und{x_{\alpha}^{p+r}z^q}v\ot\und{z^t})$ by induction. Hence, it remains to prove the required equality for $r=0$.
If $t=0$, then everything is clear. If $t>0$, then we have
\begin{multline*}
\mu_RSt_L(\und{x_{\alpha}^pz^q}v\ot\und{z^t})=\mu_R(t_L-St_Ld_Rt_L)(\und{x_{\alpha}^pz^q}v\ot\und{z^t})\\
=(-1)^{q+1}\mu_RSt_L\left(\sum\limits_{w\in G,{\gamma}\in\{0,1\}}(-1)^{jl_w+{\gamma}}\und{x_{\alpha}^{p+1}z^q}\frac{v}{x_{\alpha}}w^*\und{x_{\gamma}z^{t-1}} \frac{w}{x_{\gamma}}\right).
\end{multline*}
It is easy to see that $v=w^*=x_{\alpha}$ if $\frac{v}{x_{\alpha}}w^*=x_{\alpha}$,  and $w=\frac{v}{x_{\alpha}}$ if $\frac{v}{x_{\alpha}}w^*=1^*$. Since in both cases $ \frac{w}{x_{\alpha}}=0$, we are done by the already proved equalities.
\end{proof}

\begin{lemma}\label{zeroad} $1.$ $\mu_RSt_L(\und{x_1z^q}x_1^*\ot\und{x_{\alpha}^rz^t})=(-1)^{(q+1)(r+t)+1}\und{x_{\alpha}^rz^{q+t+1}}$ for $\alpha\in\{0,1\}$, $q,r,t\ge 0$.\\
$2.$ For $\alpha\in\{0,1\}$ and integers $r$ and $t$ one has
$$\mu_RSt_L(\und{x_1z^q}1^*\ot\und{x_{\alpha}^rz^t})=\begin{cases}
(-1)^{(q+1)(r+t)+q}x_1\und{x_{\alpha}^rz^{q+1}}+(-1)^{(r+1)q}\und{x_1^{r+2}z^q}x_1^*,&\mbox{if $t=0$ and $\alpha=1$},\\
(-1)^{(q+1)(r+t)+q}x_1\und{x_{\alpha}^rz^{q+t+1}},&\mbox{otherwise}.
\end{cases}
$$
\end{lemma}
\begin{proof} 1. Follows directly from Corollary \ref{zz}.\\
2. Can be proved analogously to Lemma \ref{zero} using Corollary \ref{zz}.
\end{proof}

\begin{lemma} If $v\in G$, ${\alpha}\in \{0,1\}$, $i>0$ and $j\ge t$, then
\begin{multline*}
Tr(\mu_RSt_L)\si_{P,P}(v\ot B_{t,j-t})=Tr(\mu_RSt_L)\si_{P,P}(v\ot C_{t,i,j-t,\alpha})\\
=Tr(\mu_RSt_L)\si_{P,P}(v\ot E_{t,i,j-t,\alpha})=0.
\end{multline*}
\end{lemma}
\begin{proof} The equality $\mu_RSt_L\si_{P,P}(v\ot B_{t,j-t})=\mu_RSt_L\si_{P,P}(v\ot E_{t,i,j-t,\alpha})=0$ follows directly from Lemma \ref{zero}. Let now prove that $\mu_RSt_L\si_{P,P}(v\ot C_{t,i,j-t,\alpha})=0$.
If $t>0$, then the required equality follows directly from Lemma \ref{zero} again. Suppose that $t=0$. Then we have
$$
Tr(\mu_RSt_L)\si_{P,P}(v\ot C_{t,i,j-t,\alpha})
=(-1)^{j+{\alpha}}\ot \sum\limits_{w\in G,\beta\in\{0,1\}}\mu_RSt_L\left(\und{x_{\alpha}^{i+1}z^{j-1}}v\frac{w^*}{x_{\alpha}}\ot \und{x_{\beta}}\frac{w}{x_{\beta}}\right).
$$
By Lemma \ref{zero} the expression $\mu_RSt_L\left(\und{x_{\alpha}^{i+1}z^{j-1}}v\frac{w^*}{x_{\alpha}}\ot \und{x_{\beta}}\frac{w}{x_{\beta}}\right)$ can be nonzero only in the case where $\alpha=\beta$ and $v\frac{w^*}{x_{\alpha}}\in\{x_{\alpha},1^*\}$.
One can show that $v\frac{w^*}{x_{\alpha}}=x_{\alpha}$ only for $v=w^*=x_{\alpha}$. Since in this case $\frac{w}{x_{\alpha}}=0$, it remains to consider the case $v\frac{w^*}{x_{\alpha}}=1^*$. In this case we have $v=wx_{\alpha}$ and so
$w=\frac{x_{\alpha}^*}{v^*}$. By Lemma \ref{zero} we have 
$$\mu_RSt_L\left(\und{x_{\alpha}^{i+1}z^{j-1}}1^*\ot \und{x_{\alpha}}\frac{x_{\alpha}^*}{v^*x_{\alpha}}\right)=(-1)^i\und{x_{\alpha}^{i+3}z^{j-1}}x_{\alpha}^*\frac{x_{\alpha}^*}{v^*x_{\alpha}}.$$
Since $\frac{x_{\alpha}^*}{v^*x_{\alpha}}\not\in\{1,x_{\alpha}\}$, we have $\mu_RSt_L\big(\und{x_{\alpha}^{i+1}z^{j-1}}v\frac{w^*}{x_{\alpha}}\ot \und{x_{\beta}}\frac{w}{x_{\beta}}\big)=0$ for all $w\in G$ and $\beta\in\{0,1\}$. Hence, we are done.
\end{proof}

\begin{lemma} If $v\in G$, ${\alpha}\in \{0,1\}$, $i>0$ and $j> t$, then
\begin{multline*}
Tr(\mu_RSt_L)\si_{P,P}(v\ot D_{t,i,j-t,\alpha})\\=\begin{cases}
0,&\mbox{if $\alpha=1$ and $v\in\{x_1,1^*\}$};\\
(-1)^{(j+1)(i+t+1)+(j-t)l_v}\frac{v}{x_{\alpha}}\ot \und{x_{\alpha}^{i+1}z^j},&\mbox{otherwise}.\\
\end{cases}
\end{multline*}
\end{lemma}
\begin{proof} Let us introduce the notation $a_{w,\beta}:=\mu_RSt_L\big(\und{x_{\beta}z^{j-t-1}}\frac{w}{x_{\beta}}v\ot\und{x_{\alpha}^{i+1}z^t}\frac{w^*}{x_{\alpha}}\big)$. Then
$$
Tr(\mu_RSt_L)\si_{P,P}(v\ot D_{t,i,j-t,\alpha})=(-1)^{(i+1)(t+1)}\ot \sum\limits_{w\in G,{\beta}\in\{0,1\}}(-1)^{(j-t)l_w+{\beta}}a_{w,\beta}.
$$
It follows from Lemmas \ref{zero} and \ref{zeroad} that if $a_{w,\beta}\not=0$, then
$\frac{w}{x_{\beta}}v\in\{x_{\beta}, x_1^*, 1^*\}$. Let us consider each of the mentioned values separately.

1. If $\frac{w}{x_{\beta}}v=x_{\beta}$, then one can show that $v=w=x_{\beta}$. Since in this case $\frac{w^*}{x_{\beta}}=0$ and $a_{w,\beta}$ can be nonzero only for $\beta=\alpha$, we obtain $a_{w,\beta}=0$.

2. If $\frac{w}{x_{\beta}}v=x_1^*$, then $a_{w,\beta}$ can be nonzero only for $\beta=1$.   One can show that $w=v^*$ in this case. On the other hand,
$
a_{v^*,1}=(-1)^{(j-t)(i+t+1)+1}\und{x_{\alpha}^{i+1}z^j}\frac{v}{x_{\alpha}}
$ if $vx_1\not=0$, and $a_{v^*,1}=0$ if $vx_1=0$.

3. If $\frac{w}{x_{\beta}}v=1^*$, then one can show that $w=x_{\beta}v^*$. One can show that $x_{\beta}^*\frac{x_{\beta}^*}{v^*x_{\beta}}=0$. Then it follows from Lemmas \ref{zero} and \ref{zeroad} that
$$
a_{x_{\beta}v^*,\beta}=\begin{cases}(-1)^{(j-t)(i+t)+1}x_1\und{x_{\alpha}^{i+1}z^j}\frac{x_1^*}{v^*x_{\alpha}},&\mbox{if $\beta=1$};\\
0,&\mbox{if $\beta=0$}.
\end{cases}
$$
Note that if $vx_1\not=0$, then $\frac{x_1^*}{v^*x_{\alpha}}=0$ and $a_{x_{\beta}v^*,1}=0$.
If $vx_1=0$, then one can show that $\frac{x_1^*}{v^*x_{\alpha}}x_1=\frac{v}{x_{\alpha}}$ except the case where $\alpha=1$ and $v\in\{x_1,1^*\}$.

Putting all the obtained equalities together we obtain the statement of the lemma.
\end{proof}

\begin{lemma}  If $v\in G$, $\alpha\in\{0,1\}$, and $i>0$, $j$ and $t$ are some integers, then
\begin{multline*}
Tr(\mu_RSt_L)\si_{P,P}\Big(\sum\limits_{r=0}^i(-1)^{rt}\und{x_{\alpha}^rz^t}\ot\und{x_{\alpha}^{i-r}z^{j-t}}\Big)\\
=(-1)^{(j+\alpha+1)(i+\alpha+1)+(j+\alpha)l_v+t(j+l_v+1)}\frac{v}{x_{\alpha}^*}\ot \und{x_{\alpha}^{i-1}z^{j+1}}+V_{i,j,t,\alpha,v}\ot \und{x_{\alpha}^{i+1}z^j},
\end{multline*}
where
$$
V_{i,j,t,\alpha,v}=\begin{cases}
0,&\mbox{if $t>0$, $\alpha=0$ and $v\in\{x_0,1^*\}$};\\
(-1)^{j(i+l_v)+1+t(j+l_v+1)}\frac{x_{\alpha}^*}{v^*},&\mbox{if $t>0$ and either $\alpha=1$ or $v\not\in\{x_0,1^*\}$};\\
(-1)^{i+j+1}\frac{v}{x_{\alpha}}+(-1)^{j(i+l_v)+1}\frac{x_{\alpha}^*}{v^*},&\mbox{if $t=0$ and $v\not\in\{x_{\alpha},1^*\}$};\\
\sum\limits_{r=1}^{i+1}(-1)^{r(i+j)+1}\frac{v}{x_{\alpha}},&\mbox{if $t=0$ and either $v=x_{\alpha}$ or  $\alpha=0$, $v=1^*$};\\
\Big(\sum\limits_{r=1}^{i}(-1)^{r(i+j)+1}+(-1)^{ij+1}\Big)\frac{v}{x_{\alpha}},&\mbox{if $t=0$, $\alpha=1$ and $v=1^*$}.
\end{cases}
$$
\end{lemma}
\begin{proof}
Using Lemmas \ref{zero} and \ref{zeroad} one can show that if $r<i$, then $$Tr(\mu_RSt_L)\si_{P,P}(v\ot\und{x_{\alpha}^rz^t}\ot\und{x_{\alpha}^{i-r}z^{j-t}})=X_{i,r,j,t,\alpha,v}+Y_{i,r,j,t,\alpha,v},$$ where
$$
X_{i,r,j,t,\alpha,v}=\begin{cases}
(-1)^{(r+1)(i+j)+1}\frac{v}{x_{\alpha}}\ot \und{x_{\alpha}^{i+1}z^j},&\mbox{if $t=0$ and either $v\in\{x_{\alpha},1^*\}$ or $r=0$};\\
0,&\mbox{otherwise},
\end{cases}
$$
and
$$
Y_{i,r,j,t,\alpha,v}=\begin{cases}
(-1)^{ji+(j+1)l_v+t(i+j+l_v)}\frac{v}{x_1^*}\ot \und{x_{1}^{i-1}z^{j+1}},&\mbox{if $r=i-1$ and $\alpha=1$};\\
0,&\mbox{otherwise}.
\end{cases}
$$

Let now calculate
\begin{multline*}
Tr(\mu_RSt_L)\si_{P,P}(v\ot\und{x_{\alpha}^iz^t}\ot\und{z^{j-t}})=1\ot \mu_R(t_L-\mu_RSt_Ld_Rt_L)(\und{z^{j-t}}v\ot\und{x_{\alpha}^{i}z^{t}})\\
=1\ot \mu_RSt_L\big(t_P(\und{z^{j-t}}v)\ot (x_{\alpha}\und{x_{\alpha}^{i-1}z^t}+(-1)^{i+t+\alpha}x_{\alpha}^*\und{x_{\alpha}^{i+1}z^{t-1}})\big).
\end{multline*}

One can show using Lemmas \ref{zero} and \ref{zeroad} that 
$$\mu_RSt_L(t_P(\und{z^{j-t}}v)x_{\alpha}\ot \und{x_{\alpha}^{i-1}z^t})=(-1)^{(j+1)(i+1)+jl_v+t(i+j+l_v+1)}(1-\alpha)\frac{v}{x_{\alpha}^*}\und{x_{\alpha}^{i-1}z^{j+1}}+Z_{i,r,j,t,\alpha,v},$$
where
$$
Z_{i,r,j,t,\alpha,v}=\begin{cases}
(-1)^{j(i+l_v)+1}\frac{x_{\alpha}^*}{v^*}\und{x_{\alpha}^{i+1}z^j},&\mbox{if $t=0$ and either $\alpha=1$ or $v\not=1^*$};\\
(-1)^{j(i+1)+1}\und{x_{0}^{i+1}z^j}x_0^*,&\mbox{if $t=0$, $\alpha=0$ and $v=1^*$};\\
0,&\mbox{otherwise}.
\end{cases}
$$
If $t>0$, then using the same lemmas one can also show that
$$
 \mu_RSt_L(t_P(\und{z^{j-t}}v)x_{\alpha}^*\ot \und{x_{\alpha}^{i+1}z^{t-1}}))=(-1)^{(i+1)(t+1)+\alpha+1}V_{i,j,t,\alpha,v}\ot \und{x_{\alpha}^{i+1}z^j}.$$

Putting all the obtained equalities together, we obtain the statement of the lemma.
\end{proof}

\begin{lemma} If $v\in G$, $\alpha\in\{0,1\}$, and $t<j$ are integers, then
\begin{multline*}
Tr(\mu_RSt_L)\si_{P,P}(v\ot A_{t,j-t})\\
=\begin{cases}
0,&\mbox{if $v\in\{1,x_1\}$};\\
(-1)^{jt+1}(x_1x_0)^i\ot\und{x_0z^j},&\mbox{if $v=x_0(x_1x_0)^i$, $0\le i\le k-1$};\\
((-1)^{jt+1}+(-1)^{j(t+1)+1})(x_0x_1)^i\ot\und{x_0z^j},&\mbox{if $v=x_1(x_0x_1)^i$, $1\le i\le k-1$};\\
(-1)^{(j+1)(t+1)}(k-i+1)x_1(x_0x_1)^{i-1}\ot\und{x_0z^j}\\
+(-1)^{(j+1)t+1}(k-i)x_0(x_1x_0)^{i-1}\ot\und{x_1z^j},&\mbox{if $v=(x_0x_1)^i$, $1\le i\le k$};\\
(-1)^{(j+1)(t+1)}(k-i)x_0(x_1x_0)^{i-1}\ot\und{x_1z^j}\\
+(-1)^{(j+1)t+1}(k-i)x_1(x_0x_1)^{i-1}\ot\und{x_0z^j},&\mbox{if $v=(x_1x_0)^i$, $1\le i\le k-1$}.
\end{cases}
\end{multline*}
\end{lemma}
\begin{proof} Let us introduce the notation $a_{u,w,\alpha,\beta}:=\mu_RSt_L\big(\und{x_{\beta}z^{j-t-1}}\frac{u}{x_{\beta}}v\frac{w^*}{u}\ot\und{x_{\alpha}z^t}\frac{w}{x_{\alpha}}\big)$. Then
$$
Tr(\mu_RSt_L)\si_{P,P}(v\ot  A_{t,j-t})=-1\ot \sum\limits_{\tiny\begin{array}{c}u,w\in G,\\ \alpha,\beta\in\{0,1\}\end{array}}(-1)^{jl_u+t(l_w+1)+{\beta}}a_{u,w,\alpha,\beta}.
$$
It follows from Lemmas \ref{zero} and \ref{zeroad} that if $a_{u,w,\alpha,\beta}\not=0$, then
$\frac{u}{x_{\beta}}v\frac{w^*}{u}\in\{x_{\beta}, x_1^*, 1^*\}$. Let us consider each of the mentioned values separately.

1. If $\frac{u}{x_{\beta}}v\frac{w^*}{u}=x_{\beta}$, then one can show that $v=u=w^*=x_{\beta}$. Since in this case $\frac{w}{x_{\beta}}=0$ and $a_{u,w,\alpha,\beta}$ can be nonzero only for $\beta=\alpha$, we obtain $a_{u,w,\alpha,\beta}=0$.

2. If $\frac{u}{x_{\beta}}v\frac{w^*}{u}=x_1^*$, then $a_{u,w,\alpha,\beta}$ can be nonzero only for $\beta=1$. Then we have $uv=wu\not=0$. Since $\frac{u}{x_{1}}\not=0$, we have $wx_1\not=0$ and $\frac{wu}{x_1}\not=0$.
Suppose that $2\nmid l_v=l_w$. Then $w=x_0(x_1x_0)^i$ for some $0\le i\le k-1$ and $u=w^*=v^*$. In this case we have $a_{v^*,v,0,1}=(-1)^{j(t+1)+1}\und{x_0z^j}\frac{v}{x_0}$, and $a_{v^*,v,1,1}=0$.

Let now $2\mid l_v=l_w$. Then $w=(x_1x_0)^i$ for some $1\le i\le k-1$ and $\frac{w}{x_{\alpha}}=0$ for $\alpha=0$. If $v=(x_1x_0)^i$, then $u=(x_1x_0)^j$ for some $1\le j\le k-i$ and we have $a_{u,w,1,1}=(-1)^{j(t+1)+1}\und{x_1z^j}\frac{v}{x_1}$.
 If $v=(x_0x_1)^i$, then $u=x_1(x_0x_1)^j$ for some $0\le j\le k-i-1$, and we have $a_{u,w,1,1}=(-1)^{j(t+1)+1}\und{x_1z^j}\frac{x_1^*}{v^*}$.

3. If $\frac{u}{x_{\beta}}v\frac{w^*}{u}=1^*$, then $wu=\frac{u}{x_{\beta}}v\not=0$. In this case we have $wx_{\beta}\not=0$. If $2\nmid l_v=l_w+1$, then we have $w=(x_{\beta}x_{\beta+1})^i$ for some $1\le i\le k-1$. Then we have either $u=x_{\beta}$, $v=x_{\beta}(x_{\beta+1}x_{\beta})^i$ or $u=(x_{\beta}x_{\beta+1})^{k-i}$, $v=x_{\beta}(x_{\beta+1}x_{\beta})^i$. We have $\frac{w}{x_{\alpha}}\not=0$ only for $\alpha=\beta$ in this case.
Note also that $\frac{1^*}{x_{\beta}}\frac{w}{x_{\beta}}=0$. Consequently, $a_{u,w,\beta,\beta}$ can be nonzero only for $\beta=1$. Hence, we have to consider only the case where $w=(x_1x_0)^i$, $v=x_1(x_0x_1)^i$. Now we have $a_{x_{1},(x_{1}x_{0})^i,1,1}=a_{(x_{1}x_{0})^{k-i},(x_{1}x_{0})^i,1,1}=(-1)^{(j+1)t+1}x_1\und{x_1z^j}x_0(x_1x_0)^{i-1}$.

Let now $2\mid l_v=l_w+1$. Then we have $w=x_{\beta+1}(x_{\beta}x_{\beta+1})^{i-1}$ for some $1\le i\le k$. In this case $\frac{w}{x_{\alpha}}\not=0$ only for $\alpha=\beta+1$ and $a_{u,w,\alpha,\beta}$ can be nonzero only for $\beta=1$, $\alpha=0$.
If $v=(x_1x_0)^{i}$ for some $i<k$, then $u=(x_1x_0)^j$ for some $1\le j\le k-i$ and we have $a_{u,w,0,1}=(-1)^{(j+1)t+1}x_1\und{x_0z^j}(x_1x_0)^{i-1}$.
 If $v=(x_0x_1)^i$ for some $i\le k$, then $u=x_1(x_0x_1)^j$ for some $0\le j\le k-i$ and we have $a_{u,w,0,1}=(-1)^{(j+1)t+1}x_1\und{x_0z^j}(x_1x_0)^{i-1}$.

Putting all the obtained equalities together we obtain the statement of the lemma.
\end{proof}

\begin{lemma}\label{last} If $v\in G$, $\alpha\in\{0,1\}$, and $j$ and $t$ are some integers, then
\begin{multline*}
Tr(\mu_RSt_L)\si_{P,P}(v\ot \und{z^t}\ot \und{z^{j-t}})\\
=\begin{cases}
0,&\mbox{if $v=x_{0}$, $t>0$};\\
(-1)^{jt+j+1}\ot\und{x_{\alpha}z^j},&\mbox{if $v=x_{\alpha}$ and either $t=0$ or $\alpha=1$};\\
(-1)^{jt+j+1}(x_{0}x_{1})^i\ot\und{x_{0}z^j},&\mbox{if $v=x_{0}(x_{1}x_{0})^i$, $1\le i\le k-1$ and $t>0$};\\
(-1)^{jt+j+1}\big((x_{\alpha}x_{\alpha+1})^i+(x_{\alpha+1}x_{\alpha})^i\big)\ot\und{x_{\alpha}z^j},&\mbox{if $v=x_{\alpha}(x_{\alpha+1}x_{\alpha})^i$, $1\le i\le k-1$}\\
&\mbox{and either $t=0$ or $\alpha=1$};\\
(-1)^{(j+1)(t+1)}(i-1)x_1(x_0x_1)^{i-1}\ot\und{x_0z^j}\\
+(-1)^{(j+1)t+1}ix_0(x_1x_0)^{i-1}\ot\und{x_1z^j},&\mbox{if $v=(x_0x_1)^i$, $1\le i\le k$, $t>0$};\\
i\big((-1)^{j+1}x_1(x_0x_1)^{i-1}\ot\und{x_0z^j}\\
-x_0(x_1x_0)^{i-1}\ot\und{x_1z^j}\big),&\mbox{if $v=(x_0x_1)^i$, $1\le i\le k$, $t=0$};\\
i\big((-1)^{(j+1)(t+1)}x_0(x_1x_0)^{i-1}\ot\und{x_1z^j}\\
+(-1)^{(j+1)t+1}x_1(x_0x_1)^{i-1}\ot\und{x_0z^j}\big),&\mbox{if $v=(x_1x_0)^i$, $0\le i\le k-1$}.\\
\end{cases}
\end{multline*}
\end{lemma}
\begin{proof} The case $t=0$ is clear. Assume now that $t>0$. Let us introduce the notation $a_{u,w,\alpha,\beta}:=\mu_RSt_L\big(\frac{u^*}{v^*}\und{x_{\alpha}z^{j-t}}\frac{u}{x_{\alpha}}w^*\ot\und{x_{\beta}z^{t-1}}\frac{w}{x_{\beta}}\big)$.
Then
\begin{multline*}
Tr(\mu_RSt_L)\si_{P,P}(v\ot \und{z^t}\ot \und{z^{j-t}})\\=
\begin{cases}
1\ot \sum\limits_{\tiny\begin{array}{c}u,w\in G\\ {\alpha,\beta}\in\{0,1\}\end{array}}(-1)^{j(l_u+l_v+1)+t(l_u+l_w+l_v+1)+\beta+1}a_{u,w,\alpha,\beta},&
\mbox{if $v\not=1^*$};\\
1\ot\sum\limits_{\tiny\begin{array}{c}u,w\in G\\ {\beta}\in\{0,1\}\end{array}}(-1)^{j(l_u+1)+t(l_u+l_w+1)+\beta+1}a_{u,w,l_u,\beta},&
\mbox{if $v=1^*$};\\
\end{cases}
\end{multline*}
It follows from Lemmas \ref{zero} and \ref{zeroad} that if $a_{u,w,\alpha,\beta}\not=0$, then
$\frac{u}{x_{\alpha}}w^*\in\{x_{\alpha}, x_1^*, 1^*\}$. Let us consider each of the mentioned values separately.

1. If $\frac{u}{x_{\alpha}}w^*=x_{\alpha}$, then one can show that $u=w^*=x_{\alpha}$. Since in this case $\frac{w}{x_{\alpha}}=0$ and $a_{u,w,\alpha,\beta}$ can be nonzero only for $\beta=\alpha$, we obtain $a_{u,w,\alpha,\beta}=0$.

2. If $\frac{u}{x_{\alpha}}w^*=x_1^*$, then $a_{u,w,\alpha,\beta}$ can be nonzero only for $\alpha=1$. In this case we have $u=w$, $\frac{w}{x_1}\not=0$. If $\beta\not=1$, then $\frac{w}{x_{\beta}}=0$.
Now we have 
$1\ot a_{u,u,1,1}=(-1)^{jt+1}\frac{u}{x_1}\frac{u^*}{v^*}\ot \und{x_1z^j}\in Tr(P)$. If $v=x_0(x_1x_0)^i$ for $0\le i\le k-1$, then $\frac{u}{x_1}\frac{u^*}{v^*}=0$ for all $u\in G$.
If $v=x_1(x_0x_1)^i$ for $0\le i\le k-1$, then
$\frac{u}{x_1}\frac{u^*}{v^*}=(x_1x_0)^i$ for $u=x_1$, $\frac{u}{x_1}\frac{u^*}{v^*}=(x_0x_1)^i$ for $u=v$, and $\frac{u}{x_1}\frac{u^*}{v^*}=0$ for all other $u\in G$.
If $v=(x_0x_1)^i$ for $1\le i\le k$, then
$\frac{u}{x_1}\frac{u^*}{v^*}=x_0(x_1x_0)^{i-1}$ for $u=x_1(x_0x_1)^j$, $0\le j\le i-1$, and $\frac{u}{x_1}\frac{u^*}{v^*}=0$ for all other $u\in G$, except the case $i=k$, $u=(x_1x_0)^j$, $1\le j\le k$ that does not occur.
If $v=(x_1x_0)^i$ for $0\le i\le k-1$, then
$\frac{u}{x_1}\frac{u^*}{v^*}=x_0(x_1x_0)^{i-1}$ for $u=(x_1x_0)^j$, $1\le j\le i$, and $\frac{u}{x_1}\frac{u^*}{v^*}=0$ for all other $u\in G$.

3. If $\frac{u}{x_{\alpha}}w^*=1^*$, then $w=\frac{u}{x_{\alpha}}$. If $\beta=\alpha$, then $\frac{w}{x_{\beta}}=0$. Hence, $a_{u,\frac{u}{x_{\alpha}},\alpha,\beta}$ can be nonzero only for $\alpha=1$, $\beta=0$.
Now we have 
$1\ot a_{u,\frac{u}{x_{1}},1,0}=(-1)^{jt+j+t}\frac{u}{x_1x_0}\frac{u^*}{v^*}x_1\ot \und{x_1z^j}\in Tr(P)$.
If $v=x_0(x_1x_0)^i$ for $0\le i\le k-1$, then
$\frac{u}{x_1x_0}\frac{u^*}{v^*}x_1=(x_0x_1)^i$ for $u=x_1x_0$, $i>0$, and $\frac{u}{x_1x_0}\frac{u^*}{v^*}x_1=0$ for all other $u\in G$.
If $v=x_1(x_0x_1)^i$ for $0\le i\le k-1$, then $\frac{u}{x_1x_0}\frac{u^*}{v^*}x_1=0$ for all $u\in G$.
If $v=(x_0x_1)^i$ for $1\le i\le k$, then
$\frac{u}{x_1x_0}\frac{u^*}{v^*}x_1=x_1(x_0x_1)^{i-1}$ for $u=x_1(x_0x_1)^j$, $1\le j\le i-1$, and $\frac{u}{x_1x_0}\frac{u^*}{v^*}x_1=0$ for all other $u\in G$, except the case $i=k$, $u=(x_1x_0)^j$, $1\le j\le k$ that does not occur.
If $v=(x_1x_0)^i$ for $0\le i\le k-1$, then
$\frac{u}{x_1x_0}\frac{u^*}{v^*}x_1=x_1(x_0x_1)^{i-1}$ for $u=(x_1x_0)^j$, $1\le j\le i$, and $\frac{u}{x_1x_0}\frac{u^*}{v^*}x_1=0$ for all other $u\in G$.

Putting all the obtained equalities together we obtain the statement of the lemma.
\end{proof}

Let us define the map $\BB_P:Tr(P)\rightarrow Tr(P[1])$ by the following equalities:
\begin{multline*}
\shoveright{\BB_P(1\ot \und{a})=0\,\,\,(a\in B);}
\end{multline*}
\begin{multline*}
\shoveright{\BB_P\big(x_{\alpha}(x_{\alpha+1}x_{\alpha})^i\ot \und{z^j}\big)=
\begin{cases}
0,&\mbox{if $0\le i\le k-1$, $2\nmid j$};\\
1\ot\und{x_{\alpha}z^j},&\mbox{if $i=0$, $2\mid j$};\\
((x_{\alpha}x_{\alpha+1})^i+(x_{\alpha+1}x_{\alpha})^i)\ot\und{x_{\alpha}z^j},&\mbox{if $1\le i\le k-1$, $2\mid j$};\\
\end{cases}}
\end{multline*}
\begin{multline*}
\BB_P\big((x_{\alpha}x_{\alpha+1})^i\ot \und{z^j}\big)=(jk+i)((-1)^{j}x_{\alpha+1}(x_{\alpha}x_{\alpha+1})^{i-1}\ot\und{x_{\alpha}z^j}+x_{\alpha}(x_{\alpha+1}x_{\alpha})^{i-1}\ot\und{x_{\alpha+1}z^j})\\
(1\le i\le k-1);
\end{multline*}
\begin{multline*}
\shoveright{\BB_P\big(1^*\ot \und{z^j}\big)=(j+1)k((-1)^{j}x_0^*\ot\und{x_{0}z^j}+x_1^*\ot\und{x_{1}z^j});}
\end{multline*}
\begin{multline*}
\shoveright{\BB_P\big(x_{\alpha+1}(x_{\alpha}x_{\alpha+1})^i\ot \und{x_{\alpha}^pz^j}\big)=0}(0\le i\le k-2);
\end{multline*}
\begin{multline*}
\shoveright{\BB_P\big(x_{\alpha}^*\ot \und{x_{\alpha}^pz^j}\big)=
\begin{cases}
(-1)^{\alpha+1}\ot \und{x_{\alpha}^{p-1}z^{j+1}},&\mbox{if $2\mid p$, $2\mid j$};\\
0,&\mbox{if $2\mid p$, $2\nmid j$};\\
(-1)^{(j+1)(\alpha+1)}(j+1)\ot\und{x_{\alpha}^{p-1}z^{j+1}},&\mbox{if  $2\nmid p$};\\
\end{cases}}
\end{multline*}
\begin{multline*}
\shoveright{\BB_P\big(x_{\alpha}\ot \und{x_{\alpha}^pz^j}\big)=
\begin{cases}
(p+1)\ot \und{x_{\alpha}^{p+1}z^j},&\mbox{if $2\mid p$, $2\mid j$;}\\
0,&\mbox{if $2\mid p$, $2\nmid j$;}\\
(-1)^{\alpha+1}j\ot \und{x_{\alpha}^{p+1}z^j},&\mbox{if $2\nmid p$, $2\mid j$;}\\
(j+p+1)\ot \und{x_{\alpha}^{p+1}z^j},&\mbox{if $2\nmid p$, $2\nmid j$;}\\
\end{cases}}
\end{multline*}
\begin{multline*}
\BB_P\big(x_{\alpha}(x_{\alpha+1}x_{\alpha})^i\ot \und{x_{\alpha}^pz^j}\big)
=
\begin{cases}
((x_{\alpha+1}x_{\alpha})^i+(x_{\alpha}x_{\alpha+1})^i)\ot \und{x_{\alpha}^{p+1}z^j},&\mbox{if $2\mid p$, $2\mid j$;}\\
0,&\mbox{if $2\mid p$, $2\nmid j$;}\\
(j+1)((-1)^{j+1}(x_{\alpha+1}x_{\alpha})^i+(x_{\alpha}x_{\alpha+1})^i)\ot \und{x_{\alpha}^{p+1}z^j},&\mbox{if $2\nmid p$;}\\
\end{cases}\\
(1\le i\le k-1);
\end{multline*}

\begin{multline*}
\BB_P\big((x_{\alpha+1}x_{\alpha})^i\ot \und{x_{\alpha}^pz^j}\big)
=
\begin{cases}
(j+1)x_{\alpha+1}(x_{\alpha}x_{\alpha+1})^{i-1}\ot \und{x_{\alpha}^{p+1}z^j},&\mbox{if $2\mid p$};\\
x_{\alpha+1}(x_{\alpha}x_{\alpha+1})^{i-1}\ot \und{x_{\alpha}^{p+1}z^j},&\mbox{if $2\nmid p$, $2\mid j$};\\
0,&\mbox{if $2\nmid p$, $2\nmid j$};\\
\end{cases}\\
(1\le i\le k-1);
\end{multline*}

\begin{multline*}
\BB_P\big((x_{\alpha}x_{\alpha+1})^i\ot \und{x_{\alpha}^pz^j}\big)
=
\begin{cases}
(-1)^{j}(j+1)x_{\alpha+1}(x_{\alpha}x_{\alpha+1})^{i-1}\ot \und{x_{\alpha}^{p+1}z^j},&\mbox{if $2\mid p$};\\
-x_{\alpha+1}(x_{\alpha}x_{\alpha+1})^{i-1}\ot \und{x_{\alpha}^{p+1}z^j},&\mbox{if $2\nmid p$, $2\mid j$};\\
0,&\mbox{if $2\nmid p$, $2\nmid j$};\\
\end{cases}\\
(1\le i\le k-1);
\end{multline*}

\begin{multline*}
\shoveright{\BB_P\big(1^*\ot \und{x_{\alpha}^pz^j}\big)=
\begin{cases}
(-1)^{\alpha}(j+1)x_{\alpha}\ot\und{x_{\alpha}^{p-1}z^{j+1}}+
(j+p+1)x_{\alpha}^*\ot \und{x_{\alpha}^{p+1}z^j},&\mbox{if $2\mid p$, $2\mid j$;}\\
(j+1)(x_{\alpha}\ot\und{x_{\alpha}^{p-1}z^{j+1}}+(-1)^{\alpha}x_{\alpha}^*\ot \und{x_{\alpha}^{p+1}z^j}),&\mbox{if  $2\mid p$, $2\nmid j$;}\\
-x_{\alpha}\ot\und{x_{\alpha}^{p-1}z^{j+1}}+px_{\alpha}^*\ot \und{x_{\alpha}^{p+1}z^j},&\mbox{if $2\nmid p$, $2\mid j$;}\\
0,&\mbox{if $2\nmid p$, $2\nmid j$.}
\end{cases}}
\end{multline*}

\begin{theorem}
The map $\BB_P$ induces the Connes' differential on $\HH_*(A)$.
\end{theorem}
\begin{proof} Follows from Lemmas \ref{zero}--\ref{last}
\end{proof}

Note that $\Hom_{A^e}(P_n,A)\cong\Homk(B_n,A)\cong A^{\dimk B_n}=A^{n+1}$. We choose this isomorphism in the following way.
We send $f\in\Hom_{A^e}(P_n,A)$ to 
$$\begin{cases}
\sum\limits_{\tiny\begin{array}{c}p+2j=n,\\ p>0,\alpha\in\{0,1\}\end{array}}f(\und{x_{\alpha}^pz^j})e_{p+\alpha}^n,&\mbox{if $2\nmid n$},\\
f(\und{z^{\frac{n}{2}}})e_{1}^n+\sum\limits_{\tiny\begin{array}{c}p+2j=n,\\ p>0,\alpha\in\{0,1\}\end{array}}f(\und{x_{\alpha}^pz^j})e_{p+\alpha}^n,&\mbox{if $2\mid n$}.
\end{cases}$$
Here $e_i^n\in A^{n+1}$ is such an element that $\pi_j^n(e_i)=0$ for $j\not=i$ and $\pi_i^n(e_i)=1$, where $\pi_j^n:A^{n+1}\rightarrow A$ ($1\le j\le n+1$) is the canonical projection on the
$j$-th component of the direct sum. We identify $\Hom_{A^e}(P_n,A)$ and $ A^{\dimk B_n}$ by the just defined isomorphism.

Let us introduce some elements of $\Hom_{A^e}(P,A)=\bigoplus\limits_{n\ge 0}A^{\dimk B_n}$.
\begin{itemize}
\item $p_1=x_0x_1+x_1x_0$, $p_2=x_1^*$, $p_2'=x_0^*$ and $p_3=1^*$ are elements of $\Hom_{A^e}(P_0,A)=A$;
\item $u_1=(x_0,0)$, $u_1'=(0,x_1)$, $u_2=(1,0)$, $u_2'=(0,1)$ are elements of $\Hom_{A^e}(P_1,A)=A^{2}$;
\item $v=(1,0,0)$, $v_1=(x_0x_1-x_1x_0,0,0)$, $v_2=(0,1,0)$, $v_2'=(0,0,1)$, $v_3=(1^*,0,0)$ are elements of $\Hom_{A^e}(P_2,A)=A^{3}$;
\item $w_1=(x_0,0,0,0)$, $w_2=(x_0^*,0,0,0)$, $w_2'=(0,x_1^*,0,0)$ are elements of $\Hom_{A^e}(P_3,A)=A^{4}$;
\item $t=(1,0,0,0,0)$ is an elements of $\Hom_{A^e}(P_4,A)=A^{5}$.
\end{itemize}

It is proved in \cite{Gen1,Gen2} that the algebra $\HH^*(A)$ is generated by the cohomological classes of the elements from $\XX$, where
$$
\XX=\begin{cases}
\{p_1,p_2,p_2',p_3,u_1,u_1',u_2,u_2',v\},&\mbox{if $\charr\kk=2$};\\
\{p_1,p_2,p_2',u_1,u_1',v_1,v_2,v_2',v_3,w_1,w_2,w_2',t\},&\mbox{if $\charr\kk\not=2$, $\charr\kk\mid k$};\\
\{p_1,p_2,p_2',u_1,u_1',v_1,v_2,v_2',t\},&\mbox{if $\charr\kk\not=2$, $\charr\kk\nmid k$}.
\end{cases}
$$
Note that our notation is essentially the same as the notation of \cite{Gen1}, but slightly differs from the notation of \cite{Gen2}. For the simplicity we denote the cohomological class of $a\in\Hom_{A^e}(P_n,A)$ by $a$ too.

It follows from the previous section that we can define the BV differential $\DD_P:\HH^*(A)\rightarrow\HH^*(A)$ by the formula $$\DD_P(f)(a)=\sum\limits_{v\in G}\theta Tr(f)\BB_P(v\ot a)v^*$$ for $a\in P$. One can show that
\begin{multline*}
\DD_P(u_2)=\DD_P(p_2'u_2)=\DD_P(u_2^2)=\DD_P(u_2')=\DD_P(p_2u_2')=\DD_P\big((u_2')^2\big)=\DD_P(v)=\DD_P(p_1v)\\
=\DD_P(u_1v)=\DD_P(u_1'v)=\DD_P(u_2v)=\DD_P(u_2'v)=\DD_P(v^2)=\DD_P(v_1)=\DD_P(p_1v_1)=\DD_P(v_2)\\
=\DD_P(p_2'v_2)=\DD_P(v_2^2)=\DD_P(v_2')=\DD_P(p_2v_2')=\DD_P\big((v_2')^2\big)=\DD_P(w_1)=\DD_P(t)=\DD_P(p_1t)\\
=\DD_P(v_1t)=\DD_P(v_2t)=\DD_P(v_2't)=\DD_P(w_1t)=\DD_P(t^2)=0,
\end{multline*}
\begin{multline*}
\DD_P(u_1)=\DD_P(u_1')=k,\,\DD_P(p_1u_1)=(k-1)p_1,\,\DD_P(p_2u_1')=p_2,\,\DD_P(p_2'u_1)=p_2',\\
\DD_P(v_3)=u_1'-u_1\,,\DD_P(p_2v)=u_2',\,\DD_P(p_2'v)=u_2,\,\DD_P(p_3v)=u_1+u_1',\\
\DD_P(u_1u_1')=k(u_1'-u_1),\,\DD_P(u_1u_2)=ku_2,\,\DD_P(u_1'u_2')=ku_2',\,\DD_P(w_2)=v_2,\\
\DD_P(w_2')=-v_2',\,\DD_P(u_1v_1)=(2k-1)v_1,\,\DD_P(u_1v_2)=(k+2)v_2,\,\DD_P(u_1v_2')=kv_2',\\
\DD_P(u_1'v_2)=kv_2,\,\DD_P(u_1'v_2')=(k+2)v_2',\,\DD_P(v_2v_3)=3w_2,\,\DD_P(v_2'v_3)=3w_2',\\
\DD_P(u_1t)=\DD_P(u_1't)=3kt,\,\DD_P(v_2w_2)=v_2^2,\,\DD_P(v_2'w_2')=-(v_2')^2,\\
\DD_P(v_3t)=3(u_1'-u_1)t,\,\DD_P(w_2t)=3v_2t,\,\DD_P(w_2't)=3v_2't.
\end{multline*}
We use here the results of \cite{Gen1,Gen2}. In particular, we use the formulas for some products in $\HH^*(A)$ and the description of some coboundaries. Alternatively, one can use the formula $f\smile g=(f\ot g)S\eta_L$ and Lemmas \ref{zero} and \ref{zeroad} to compute products in $\HH^*(A)$. Note also that in each of the formulas above we assume that all the elements included in the formula lie in $\XX$. For example, if $v$ appears in some equality, then this equality holds for $\charr\kk=2$, but doesn't have to hold for $\charr\kk\not=2$. We also have $\DD_P(a)=0$ for all $a\in\HH^0(A)$. Now it is not hard to recover the Gerstenhaber bracket and the rest of the BV differential on the Hochschild cohomology of $A$ using relations between the generators of $\HH^*(A)$ described in \cite{Gen1,Gen2} and the graded Leibniz rule for the Gerstenhaber bracket and the cup product.

\begin{Remark}{\rm We think that the Gerstenhaber bracket and even the BV differential on $\HH^*(A)$ can be computed in some simpler way using different tricks. But the aim of this example is to show that our formulas are reasonable for a direct application.}
\end{Remark}

%Using results of \cite{Gen1,Gen2}, it is easy to see that we have defined $\DD_P$ on the classes of the elements of $\XX$ and on the elements of the form $ab$, where the $a$ and $b$ are cohomology classes of elements from $\XX$.
{\bf Acknowledgements.} The author is grateful to Sergey Ivanov, Maria Julia Redondo, Dmitry Kaledin, and especially to Sarah Witherspoon for productive discussions, helpful advises, and the attention to my work.

\end{document}